\documentclass[10pt]{article}

\usepackage{amsmath}
\usepackage{amsthm}
\usepackage{amssymb}
\usepackage{amscd}
\usepackage{amsfonts}  
\usepackage{graphicx}    
\usepackage{multicol}    
\input xypic
\xyoption{all}

\DeclareMathAlphabet\EuR{U}{eur}{m}{n}
\SetMathAlphabet\EuR{bold}{U}{eur}{b}{n}

\begin{document}


\newtheorem{theorem}{Theorem}[section]
\newtheorem{lemma}[theorem]{Lemma}
\newtheorem{proposition}[theorem]{Proposition}
\newtheorem{definition}[theorem]{Definition}
\newtheorem{example}[theorem]{Example}
\newtheorem{remark}[theorem]{Remark}
\newtheorem{corollary}[theorem]{Corollary}
\newtheorem{conjecture}[theorem]{Conjecture}
\newtheorem{problem}[theorem]{Problem}
\newtheorem{question}[theorem]{Question}
{\catcode`@=11\global\let\c@equation=\c@theorem}
\renewcommand{\theequation}{\thetheorem}

\renewcommand{\theenumi}{\alph{enumi}}
\renewcommand{\labelenumi}{(\theenumi)}

\makeatletter
\renewcommand{\@seccntformat}[1]{\csname the#1\endcsname.\hspace{1em}}
\makeatother
\newcommand{\tit}[1]{\begin{bf} \begin{center} \begin{Large}
\section{#1}
\label{sec: #1}
\end{Large}\end{center}\end{bf}
\nopagebreak}


\newcommand{\squarematrix}[4]{\left( \begin{array}{cc} #1 & #2 \\ #3 &
#4
\end{array} \right)}
\newcommand{\smallmx}[4]{\mbox{\begin{scriptsize}$\squarematrix{#1}{#2}
        {#3}{#4}$\end{scriptsize}}}

\newcommand{\comsquare}[8]{
\begin{center}
$\begin{CD}
#1 @>#2>> #3\\
@V{#4}VV @VV{#5}V\\
#6 @>>#7> #8
\end{CD}$
\end{center}}


\let\sect=\S
\newcommand{\curs}{\EuR}
\newcommand{\CHAINCOMPLEXES}{\curs{CHCOM}}
\newcommand{\GROUPOIDS}{\curs{GROUPOIDS}}
\newcommand{\PAIRS}{\curs{PAIRS}}
\newcommand{\FGINJ}{\curs{FGINJ}}
\newcommand{\FGMOD}{\curs{FGMOD}}
\newcommand{\FGFMOD}{\curs{FGFMOD}}
\newcommand{\FGPMOD}{\curs{FGPMOD}}
\newcommand{\MOD}{\curs{MOD}}
\newcommand{\Or}{\curs{Or}}
\newcommand{\FSETS}{\curs{FSETS}}
\newcommand{\SPACES}{\curs{SPACES}}
\newcommand{\SPECTRA}{\curs{SPECTRA}}
\newcommand{\Sub}{\curs{Sub}}


\newcommand{\bbC}{{\mathbb C}}
\newcommand{\bbH}{{\mathbb H}}
\newcommand{\bbI}{{\mathbb I}}
\newcommand{\bbK}{{\mathbb K}}
\newcommand{\bbKO}{\mathbb{KO}}
\newcommand{\bbN}{{\mathbb N}}
\newcommand{\bbP}{{\mathbb P}}
\newcommand{\bbQ}{{\mathbb Q}}
\newcommand{\bbR}{{\mathbb R}}
\newcommand{\bbZ}{{\mathbb Z}}

\newcommand{\calbh}{{\mathcal B}{\mathcal H}}
\newcommand{\calc}{{\mathcal C}}
\newcommand{\cald}{{\mathcal D}}
\newcommand{\cale}{{\mathcal E}}
\newcommand{\calf}{{\mathcal F}}
\newcommand{\calg}{{\mathcal G}}
\newcommand{\calh}{{\mathcal H}}
\newcommand{\calk}{{\mathcal K}}
\newcommand{\calp}{{\mathcal P}}
\newcommand{\cals}{{\mathcal S}}

\newcommand{\calall}{{\mathcal A} {\mathcal L}{\mathcal L}}
\newcommand{\calfcyc}{{\mathcal F}{\mathcal C}{\mathcal Y}}
\newcommand{\calcyc}{{\mathcal C}{\mathcal Y}{\mathcal C}}
\newcommand{\calfin}{{\mathcal F}{\mathcal I}{\mathcal N}}
\newcommand{\calko}{{\mathcal K}{\mathcal O}}
\newcommand{\calmfin}{{\mathcal M}{\mathcal F}{\mathcal I}}
\newcommand{\calvcyc}{{\mathcal V}{\mathcal C}{\mathcal Y}}
\newcommand{\caltr}{{\cal T}{\cal R}}
\newcommand{\trivial}{ \{ \! 1 \! \} }

\newcommand{\bfA}{\ensuremath{\mathbf{A}}}
\newcommand{\bfa}{\ensuremath{\mathbf{a}}}
\newcommand{\bfb}{\ensuremath{\mathbf{b}}}
\newcommand{\bfCTR}{\ensuremath{\mathbf{CTR}}}
\newcommand{\bfDtr}{\ensuremath{\mathbf{Dtr}}}
\newcommand{\bfE}{\ensuremath{\mathbf{E}}}
\newcommand{\bff}{\ensuremath{\mathbf{f}}}
\newcommand{\bfF}{\ensuremath{\mathbf{F}}}
\newcommand{\bfg}{\ensuremath{\mathbf{g}}}
\newcommand{\bfHH}{\ensuremath{\mathbf{HH}}}
\newcommand{\bfHPB}{\ensuremath{\mathbf{HPB}}}
\newcommand{\bfI}{\ensuremath{\mathbf{I}}}
\newcommand{\bfi}{\ensuremath{\mathbf{i}}}
\newcommand{\bfK}{\ensuremath{\mathbf{K}}}
\newcommand{\bfL}{\ensuremath{\mathbf{L}}}
\newcommand{\bfr}{\ensuremath{\mathbf{r}}}
\newcommand{\bfs}{\ensuremath{\mathbf{s}}}
\newcommand{\bfS}{\ensuremath{\mathbf{S}}}
\newcommand{\bft}{\ensuremath{\mathbf{t}}}
\newcommand{\bfT}{\ensuremath{\mathbf{T}}}
\newcommand{\bfTC}{\ensuremath{\mathbf{TC}}}
\newcommand{\bfU}{\ensuremath{\mathbf{U}}}
\newcommand{\bfu}{\ensuremath{\mathbf{u}}}
\newcommand{\bfv}{\ensuremath{\mathbf{v}}}
\newcommand{\bfw}{\ensuremath{\mathbf{w}}}


\newcommand{\aut}{\operatorname{aut}}
\newcommand{\Bor}{\operatorname{Bor}}
\newcommand{\ch}{\operatorname{ch}}
\newcommand{\class}{\operatorname{class}}
\newcommand{\cok}{\operatorname{coker}}
\newcommand{\cone}{\operatorname{cone}}
\newcommand{\colim}{\operatorname{colim}}
\newcommand{\COMPLEXES}{\operatorname{COMPLEXES}}
\newcommand{\con}{\operatorname{con}}
\newcommand{\conhom}{\operatorname{conhom}}
\newcommand{\cyclic}{\operatorname{cyclic}}
\newcommand{\ev}{\operatorname{ev}}
\newcommand{\Ext}{\operatorname{Ext}}
\newcommand{\func}{\operatorname{func}}
\newcommand{\Gen}{\operatorname{Gen}}
\newcommand{\Hom}{\operatorname{Hom}}
\newcommand{\hur}{\operatorname{hur}}
\newcommand{\im}{\operatorname{im}}
\newcommand{\inj}{\operatorname{inj}}
\newcommand{\id}{\operatorname{id}}
\newcommand{\infl}{\operatorname{Infl}}
\newcommand{\ind}{\operatorname{ind}}
\newcommand{\Inn}{\operatorname{Inn}}
\newcommand{\Irr}{\operatorname{Irr}}
\newcommand{\Is}{\operatorname{Is}}
\newcommand{\map}{\operatorname{map}}
\newcommand{\mor}{\operatorname{mor}}
\newcommand{\Ob}{\operatorname{Ob}}
\newcommand{\op}{\operatorname{op}}
\newcommand{\pr}{\operatorname{pr}}
\newcommand{\point}{\operatorname{pt.}}
\newcommand{\Rat}{\operatorname{Rat}}
\newcommand{\res}{\operatorname{res}}
\newcommand{\Sw}{\operatorname{Sw}}
\newcommand{\topo}{\operatorname{top}}
\newcommand{\Tor}{\operatorname{Tor}}
\newcommand{\tors}{\operatorname{tors}}
\newcommand{\Wh}{\operatorname{Wh}}

\newcommand{\pt}{\{\operatorname{pt}. \}}

\newcounter{commentcounter}
\newcommand{\comment}[1]                      
{\stepcounter{commentcounter}
{\bf Comment \arabic{commentcounter}}: {\ttfamily #1} }

\newcommand{\EFG}[2]{E_{#1}#2}               
\newcommand{\OrF}[2]{\Or_{#1}(#2)}               
\newcommand{\version}[1]{\begin{center} last edited on #1 or later\\
last compiled on \today \end{center}} 

\title{Induction Theorems and Isomorphism Conjectures for $K$- and $L$-Theory}
\author{Arthur Bartels\thanks{\noindent email:
bartelsa@math.uni-muenster.de} and
Wolfgang L\"uck\thanks{\noindent email:
lueck@math.uni-muenster.de\protect\\
www: ~http://www.math.uni-muenster.de/u/lueck/\protect\\
fax: +49 251 8338370\protect}
\\
Fachbereich Mathematik\\ Universit\"at M\"unster\\
Einsteinstr.~62\\ 48149 M\"unster\\Germany}

\maketitle

\typeout{-----------------------  Abstract  ------------------------}

\begin{abstract}
The Farrell-Jones and the Baum-Connes Conjecture say that one can compute the algebraic
$K$- and $L$-theory of the group ring and the topological  $K$-theory of
the reduced group $C^*$-algebra of a group $G$
in terms of these functors for the virtually cyclic subgroups or the finite subgroups of $G$.
By induction theory we want to reduce these families of subgroups to a smaller family,
for instance to the family of subgroups which are either finite hyperelementary or
extensions of finite hyperelementary groups
with $\bbZ$ as kernel or to the family of finite cyclic subgroups.
Roughly speaking, we extend the induction theorems of Dress for finite groups
to infinite groups.

\smallskip


\end{abstract}
\vspace{4mm}


\typeout{-----------------------  Introduction  ------------------------}
\setcounter{section}{-1}
\section{Introduction}
\label{sec:Introduction}
The goal of this paper is to reduce the family of virtually cyclic subgroups
to  a smaller subfamily in the statements of the (Fibered) Farrell-Jones Conjecture for
algebraic $K$- and $L$-theory of group rings. The strategy is to extend
the classical induction results for finite groups to infinite groups.

Let $\calfin$ and $\calvcyc$ respectively be the class of finite groups and virtually cyclic groups respectively.
Let $\calf \subseteq \calfin$ be a subclass of the class of finite groups
which is closed under isomorphism of groups and taking subgroups. Define
\begin{eqnarray*}
\calf' &\subseteq & \calvcyc
\end{eqnarray*}
to be the class of groups $V$ for which either

\begin {enumerate}
\item there exists an extension $1\to \bbZ \to V \to F \to 1$ for a group $F \in \calf$ or
\item $V \in \calf$ holds.
\end{enumerate}
With this notion we get $\calvcyc = \calfin'$.

Let $p$ be a prime. A finite group $G$ is called \emph{$p$-elementary}
if it is isomorphic to $C \times P$ for a cyclic group $C$ and a $p$-group $P$
such that the order $|C|$ is prime to $p$.
A finite group $G$ is called \emph{$p$-hyperelementary}
if it can be written as an extension
$1 \to C \to G \to P\to 1$ for a cyclic group $C$ and a $p$-group $P$
such that the order $|C|$ is prime to $p$.
A finite group $G$ is called \emph{elementary} or \emph{hyperelementary} respectively
if it is \emph{$p$-elementary} or \emph{$p$-hyperelementary} respectively
for some prime $p$.
Let $\calfcyc$ be the class of finite cyclic groups. Let $\cale_p$ and $\calh_p$ respectively be
the class of groups which are
$p$-elementary groups and $p$-hyperelementary respectively for a prime $p$. Let $\cale$ and $\calh$
respectively be the class of
groups which are elementary and hyperelementary respectively.

A family of subgroups of $G$ is a set of subgroups which is closed under conjugation and taking subgroups.
For a class $\calf$ of groups which is closed under isomorphism of groups and taking subgroups
we denote by $\calf(G)$ the family of subgroups of $G$ whose members are in $\calf$. Thus $\calfin(G)$ is
the family of finite subgroups of $G$. If $G$ is clear from the context we will often
abuse notation and write simply $\calf$ for $\calf(G)$.

The Isomorphism Conjecture of Farrell and Jones appeared in \cite[1.6]{Farrell-Jones(1993a)}.
For a survey of this conjecture see \cite{Lueck-Reich(2004g)}.

\begin{theorem}[Induction theorem for algebraic $K$-theory]
\label{the: Induction theorem for algebraic K-theory}
Let $G$ be a group and let $N$ be an integer. Then the following hold.

\begin{enumerate}

\item \label{the: Induction theorem for algebraic K-theory: calh}
The group $G$ satisfies the (Fibered) Isomorphism Conjecture (in the range $\le N$) for algebraic $K$-theory
with coefficients in $R$ for the family $\calvcyc$ if and only if
$G$ satisfies the (Fibered) Isomorphism Conjecture (in the range $\le N$) for algebraic $K$-theory
with coefficients in $R$ for the family $\calh'$.

\item \label{the: Induction theorem for algebraic K-theory: calh_p}
Let $p$ be a prime. Then $G$ satisfies the (Fibered) Isomorphism Conjecture
(in the range $\le N$) for algebraic $K$-theory
with coefficients in $R$ for the family $\calvcyc$ after applying
$ \bbZ_{(p)}\otimes_{\bbZ} -$ if and only if $G$  satisfies the (Fibered) Isomorphism Conjecture
(in the range $\le N$) for algebraic $K$-theory
with coefficients in $R$ for the family $\calh_p'$ after applying
$\bbZ_{(p)} \otimes_{\bbZ} -$.

\item \label{the: Induction theorem for algebraic K-theory: calh and Q subset R, R regular}
Suppose that $R$ is regular and $\bbQ \subseteq R$.
Then the group $G$ satisfies the  Isomorphism Conjecture (in the range $\le N$) for algebraic $K$-theory
with coefficients in $R$ for the family $\calvcyc$ if and only if
$G$ satisfies the Isomorphism Conjecture (in the range $\le N$) for algebraic $K$-theory
with coefficients in $R$ for the family $\calh$.

If we assume that $R$ is regular and $\bbC \subseteq R$, then we can
replace $\calh$ by $\cale$.

\item \label{the: Induction theorem for algebraic K-theory: calh_p and Q subset R, R regular}
Suppose that $R$ is regular and $\bbQ \subseteq R$.
Let $p$ be a prime. Then $G$ satisfies the
Isomorphism Conjecture (in the range $\le N$) for algebraic $K$-theory
with coefficients in $R$ for the family $\calvcyc$ after applying
$ \bbZ_{(p)}\otimes_{\bbZ} -$ if and only if $G$  satisfies the
Isomorphism Conjecture (in the range $\le N$) for algebraic $K$-theory
with coefficients in $R$ for the family $\calh_p$ after applying
$\bbZ_{(p)} \otimes_{\bbZ} -$.

If we assume that $R$ is regular and $\bbC \subseteq R$, then we can
replace $\calh_p$ by $\cale_p$.

\end{enumerate}
\end{theorem}

Note that in \eqref{the: Induction theorem for algebraic K-theory: calh and Q subset R, R regular} and
\eqref{the: Induction theorem for algebraic K-theory: calh_p and Q subset R, R regular} above there are no claims about
the Fibered Isomorphism Conjecture. The problem is that we do not know
whether the relevant Nil-groups for amalgamated products and
HNN-extensions with rational coefficients vanish for all groups.
If the group $\bbZ$ satisfies the Fibered Isomorphism Conjecture for algebraic $K$-theory with coefficients in
$\bbQ$ for the trivial family, then the Nil groups $NK_*(\bbQ \Gamma)$ vanishes for all groups $\Gamma$.
It seems to be not known whether these Nil groups vanish.
(Of course for the group $\bbZ$ the families $\calh$, $\calh_p$ and $\calfin$ are all trivial.)

\begin{remark}[Relative assembly maps] \label{rem:relative_assembly_maps}
\em
We have stated the main theorems in terms of the (Fibered) Isomorphism Conjectures.
The main new result we prove is that the relevant relative assembly maps are bijective.
This statement is true in general and is independent of the question whether
the Isomorphism Conjecture is true or not but implies our main results.

For instance, the assembly map for the family $\calh'$ factorizes as
$$\calh_n(\EFG{\calh'}{G}) \to \calh_n(\EFG{\calfcyc}{G}) \to \calh_n^G(\pt) = K_n(RG)$$
where $\calh^G_*$ is a certain $G$-homology theory related to algebraic $K$-theory of group rings with coefficients
in a ring $R$,  the first map is a relative assembly map which we will prove is bijective, and the second map
is the assembly map for $\calvcyc$. The Isomorphism Conjecture for algebraic $K$-theory
for $\calh'$ or $\calvcyc$ respectively says that the assembly map for $\calh'$ or $\calvcyc$
respectively is bijective for $n \in \bbZ$.
The same remark applies to the $L$-theory version below.
\em
\end{remark}

\begin{theorem}[Induction theorem for algebraic $L$-theory]
\label{the: Induction theorem for algebraic L-theory}
Let $G$ be a group. Then the following hold.

\begin{enumerate}

\item \label{the: Induction theorem for algebraic L-theory: calh_2 cup_{p not= 2} cale_p}
The group $G$ satisfies the (Fibered) Isomorphism Conjecture (in the range $\le N$) for algebraic $L$-theory
with coefficients in $R$ for the family $\calvcyc$ if and only if
$G$ satisfies the (Fibered) Isomorphism Conjecture for algebraic $L$-theory
with coefficients in $R$ for the family $\left(\calh_2 \cup \bigcup_{p\text{ prime}, p \not = 2} \cale_p\right)'$.

\item \label{the: Induction theorem for algebraic L-theory: cup_{p not= 2} cale_p}
The group $G$ satisfies the  (Fibered) Isomorphism Conjecture (in the range $\le N$) for algebraic $L$-theory
with coefficients in $R$ for the family $\calvcyc$ after applying
$ \bbZ[1/2]\otimes_{\bbZ} -$ if and only if  $G$  satisfies the
(Fibered) Isomorphism Conjecture for algebraic $L$-theory
with coefficients in $R$ for the family $\bigcup_{p\text{ prime}, p \not = 2} \cale_p$
after applying $\bbZ[1/2] \otimes_{\bbZ} -$.

\end{enumerate}
\end{theorem}

\begin{remark}[Rationalized versions] \label{rem: rationalized versions} \em
We omit the discussion of the rationalized versions, i.e.\ the versions obtained after applying
$\bbQ \otimes_{\bbZ} -$. In this case one gets more precise information as discussed in detail in
\cite[Section 8]{Lueck-Reich(2004g)}. The results presented there are based on
\cite{Bartels(2003h)}, \cite{Grunewald(2005)}, \cite{Kuku-Tang(2003)} and \cite{Lueck(2002b)}.
\em
\end{remark}

%
%

The next result is due to Mislin and Matthey
\cite{Matthey-Mislin(2004)} for the complex case.  We will give a proof
for both the complex and the real case in our framework.
It is not clear to us whether it is possible to extend the methods of \cite{Matthey-Mislin(2004)}
to the real case.

\begin{theorem}[Induction theorem for topological $K$-theory]
\label{the: Induction theorem for topological K-theory}
Let $G$ be a group.

Then the relative assembly map
$$K_n^G(\EFG{\calfcyc}{G})\to K_n^G(\EFG{\calfin}{G})$$
is bijective for all $n \in \bbZ$.

In particular $G$ satisfies the Baum-Connes Conjecture if and only if
$G$ satisfies the Baum-Connes Conjecture for the family $\calfcyc$.

The corresponding statements are also true if one replaces complex equivariant $K$-homology  $K^G_*$ by
real equivariant $K$-homology  $KO^G_*$ and uses the real reduced group
$C^*$-algebra in the Baum-Connes Conjecture.
\end{theorem}

The rational version of the Induction theorem for topological real $K$-theory has been
used by Stolz \cite[p.695]{Stolz(2002)}.
Using Artin Induction \cite[Theorem~26 on page 97]{Serre(1977)} the methods we use
to proof the Induction Theorem for algebraic K-theory
can also be used to give a simpler proof of the rational version of
Theorem~\ref{the: Induction theorem for topological K-theory}.

The paper is organized as follows\\[2mm]
\begin{tabular}{ll}
\ref{sec: Transitivity Principles}. &  Transitivity Principles
\\
\ref{sec: General Induction Theorems}. &  General Induction Theorems
\\
\ref{sec: The Swan Group as a Functor on Groupoids}. & The Swan Group as a Functor on Groupoids
\\
\ref{sec: Proof of the Main Result for Algebraic $K$-Theory}. &  Proof of the Main Result for Algebraic $K$-Theory
\\
\ref{sec: Outline of the Proof of the Main Result for Algebraic $L$-Theory}. &
Outline of the Proof of the Main Result for Algebraic $L$-Theory
\\
\ref{sec: Proof of the Main Result for Topological  $K$-Theory}. &  Proof of the Main Result for Topological $K$-Theory
\\
\ref{sec: Versions in terms of colimits}. & Versions in terms of colimits
\\
\ref{sec: On Quinn's Hyperelementary Induction Conjecture}. & On Quinn's Hyperelementary Induction Conjecture
\\
&  References
\end{tabular}

\typeout{--------------------------  Section 1  ------------------------}
\section{Transitivity Principles}
\label{sec: Transitivity Principles}

In this section we fix an equivariant homology theory  $\calh^?_*$ with values in
$\Lambda$-modules for a commutative associative ring $\Lambda$ with
unit in the sense of
\cite[Section 1]{Lueck(2002b)}. This essentially means that
we get for each group $G$ a $G$-homology theory $\calh^G_*$ which
assigns to a (not necessarily proper or cocompact) pair of
$G$-$CW$-complexes $(X,A)$ a $\bbZ$-graded
$\Lambda$-module $\calh^G_n(X,A)$, such that there exists natural
long exact sequences of pairs and $G$-homotopy invariance, excision,
and the disjoint union axiom are satisfied. Moreover, an induction
structure is required which in particular implies for a subgroup
$H \subseteq G$
and a $H$-$CW$-pair $(X,A)$ that there is a natural isomorphism
$\calh^H_n(X,A) \xrightarrow{\cong} \calh^G_n(G \times_H(X,A))$.

Recall that a family $\calf$ of subgroups of $G$ is a set of subgroups which is
closed under conjugation and taking subgroups. Examples are the families $\calfin(G)$
of finite subgroups and $\calvcyc(G)$ of virtually cyclic subgroups. Given a group homomorphism
$\phi \colon K \to G$ and a family $\calf$ of subgroups of $G$,
define the family
$\phi^*\calf$ of subgroups of $K$ by
$$\phi^*\calf ~ := \{H \subseteq K \mid \phi(H) \in \calf\}.$$
If $\phi$ is an inclusion of subgroups, we also write
$$\phi^*\calf = K \cap \calf ~ = ~  \{H \subseteq K \mid H \in \calf\}
~ =  ~ \{L \cap K \mid L \in \calf\}.$$
If $\psi \colon H \to K$ is another group homomorphism, then
\begin{eqnarray}
\psi^*(\phi^* \calf) & = & (\phi \circ \psi)^*\calf.
\label{psi^*phi^*calf= (psi circ psi)^*calf}
\end{eqnarray}

Associated to a family $\calf$ there is a $G$-$CW$-complex
$\EFG{\calf}{G}$ (unique up to $G$-homotopy equivalence)
with the property that the fixpoint sets $(\EFG{\calf}{G})^H$
are contractible for $H \in \calf$ and
empty for $H \notin \calf$. It is called the \emph{classifying space of the family $\calf$}.
For more information about these spaces we refer for instance to \cite{Lueck(2004h)}, \cite[pp.46]{Dieck(1987)}.

\begin{definition}[(Fibered) Isomorphism Conjecture]
\label{def: (Fibered) Isomorphism Conjectures for calh^?_*}
Fix an equivariant homology theory  $\calh^?_*$ with values in
$\Lambda$-modules for a commutative associative ring $\Lambda$.
A group $G$ together with a family  of subgroups $\calf$ satisfies the
\emph{Isomorphism Conjecture (in the range $\le N$)}
if the projection $\pr \colon \EFG{\calf}{G} \to \pt$
 induces an isomorphism
$$\calh^G_n(\pr) \colon \calh^G_n(\EFG{\calf}{G})  \xrightarrow{\cong} \calh^G_n(\pt)$$
 for $n \in \bbZ$  (with $n \le N)$.

The pair $(G,\calf)$  satisfies  the \emph{Fibered Isomorphism
 Conjecture (in the range $\le N$)} if for each group homomorphism
$\phi \colon K \to G$ the pair $(K,\phi^*\calf)$ satisfies the
 Isomorphism Conjecture (in the range $\le N$).
\end{definition}

Built in into the Fibered Isomorphism Conjecture is the
following obvious inheritance
property which is not true in general in the non-fibered case.

\begin{lemma} \label{lem: basic inheritance property of fibered conjecture}
Let $\phi \colon K \to G$ be a group
homomorphism and let $\calf$ a family of subgroups. If $(G,\calf)$ satisfies
the Fibered Isomorphism Conjecture
\ref{def: (Fibered) Isomorphism Conjectures for calh^?_*},
then $(K,\phi^*\calf)$ satisfies
the Fibered Isomorphism Conjecture
\ref{def: (Fibered) Isomorphism Conjectures for calh^?_*}.
\end{lemma}
\begin{proof}
If $\psi \colon L \to K$ is a group homomorphism, then
$\psi^*(\phi^*\calf) = (\phi \circ \psi)^*\calf$
from \eqref{psi^*phi^*calf= (psi circ psi)^*calf}.
\end{proof}

We will use several times the following result proved in \cite[Theorem 2.9]{Lueck-Reich(2004g)}

\begin{theorem}[Transitivity Principle for equivariant homology] \label{the: transitivity principle of bcsfinal}
Suppose $\calf \subset \calg$ are two families of subgroups of the group $G$.
Let $N$ be an integer. If for every $H \in \calg$ and every $n \leq N$ the map
induced by the projection
$$
\calh_n^H (\EFG{\calf \cap H}{H}) \to \calh_n^H (\pt)
$$
is an isomorphism, then for every $n \leq N$ the map induced by the up to $G$-homotopy
unique $G$-map $\EFG{\calf}{G} \to \EFG{\calg}{G}$
$$
\calh_n^G (\EFG{\calf}{G} ) \to \calh_n^G ( \EFG{\calg}{G} )
$$
is an isomorphism.
\end{theorem}

This implies the following transitivity principle for the (Fibered) Isomorphism Conjecture.
At the level of spectra this transitivity principle is due to
Farrell and Jones \cite[TheoremA.10]{Farrell-Jones(1993a)}.

\begin{theorem}[Transitivity Principle] \label{the: transitivity}
Suppose $\calf \subseteq \calg$ are two families of subgroups of $G$.
Assume that for every element $H \in \calg$ the group $H$ satisfies
the (Fibered) Isomorphism Conjecture for $\calf \cap H$
 (in the range $\le N$).

Then $(G,\calg)$ satisfies the (Fibered) Isomorphism Conjecture
(in the range $\le N$)
if and only if  $(G,\calf)$ satisfies the (Fibered) Isomorphism
Conjecture (in the range $\le N$).
\end{theorem}

\begin{proof}
We first prove the claim for the Fibered Isomorphism Conjecture.
Consider a group homomorphism $\phi \colon K \to G$. Then for every
subgroup $H$ of $K$ we conclude
$$(\phi|_H)^*(\calf \cap \phi(H)) ~ = ~ (\phi^*\calf) \cap H$$
from \eqref{psi^*phi^*calf= (psi circ psi)^*calf}, where $\phi|_H \colon H \to \phi(H)$ is  the
group homomorphism induced by $\phi$.
For every element $H \in \phi^*\calg$ the map
$$\calh_p^H(\EFG{(\phi|_H)^*(\calf \cap \phi(H))}{H}) =
\calh_p^H(\EFG{\phi^*\calf \cap H}{H}) \to
\calh_p^H(\pt)$$
is bijective for all $p \in \bbZ$ ( with $p \le N)$ by the assumption that
that the element $\phi(H) \in \calg$  satisfies
the Fibered Isomorphism Conjecture for $\calf \cap \phi(H)$ (in the range $\le N$).
Hence by
Theorem~\ref{the: transitivity principle of bcsfinal} applied to the inclusion
$\phi^*\calf \subseteq \phi^*\calg$  of families of subgroups of $K$
we get an isomorphism
$$\calh_p^K(\EFG{\phi^*\calf}{K}) \xrightarrow{\cong}
\calh_p^K(\EFG{\phi^*\calg}{K})$$
(for $p \le N)$. Therefore the map
$\calh_p^K(\EFG{\phi^*\calf}{K}) \to \calh_p^K(\pt)$
is bijective for all
$p \in \bbZ$ (with $p \le N)$
if and only if the map
$\calh_p^K(\EFG{\phi^*\calg}{K}) \to \calh^K_p(\pt)$
is bijective for all
$p \in \bbZ$ (with $p \le N)$.

If we specialize this argument to $\phi = \id_G$ then only the Isomorphism Conjecture for $\calf \cap H$
(in the range $\le N$) is needed and this proves the claim for the Isomorphism Conjecture.
\end{proof}

As an easy application we see that we can always enlarge the family in the Fibered Isomorphism Conjecture.
Note that this argument is not valid for the Isomorphism Conjecture.

\begin{lemma} \label{lem: enlarging the family for the Fibered Isomorphism Conjecture}
Let $G$ be a group and let $\calf \subset \calg$ be families of subgroups of $G$.
Suppose that $G$ satisfies the Fibered Isomorphism
Conjecture~\ref{def: (Fibered) Isomorphism Conjectures for calh^?_*} (in the range $n \le N$) for the family $\calf$.

Then $G$ satisfies the Fibered Isomorphism
Conjecture~\ref{def: (Fibered) Isomorphism Conjectures for calh^?_*} (in the range $n \le N$) for the family $\calg$.
\end{lemma}
\begin{proof}
We want to use
Theorem~\ref{the: transitivity}. Therefore we have to show for each subgroup $K$ in $\calg$
that it satisfies the Fibered Isomorphism Conjecture for $\calf \cap K$.
If $i \colon K \to G$ is the inclusion, then
$i^*\calf = \calf \cap K$. Now apply
Lemma~\ref{lem: basic inheritance property of fibered conjecture}.
\end{proof}

We have defined in the introduction for every class $\calf$ of finite groups
the class $\calf'$ consisting of those virtually cyclic groups $V$
for which there exists an extension $1 \to \bbZ \to V \to F \to 1$ for a group $F \in \calf$
or for which $V \in \calf$ holds.

\begin{theorem}\label{the: calfin to calvcyc}
Let $\calf$ be a class of finite groups closed under isomorphism and taking subgroups.
Suppose that every finite group $F$ satisfies the Fibered Isomorphism Conjecture
(in the range $\le N$)
with respect to the family $\calf(F) = \{H \subseteq F\mid H \in \calf\}$.
Let $G$ be a group.

Then $G$ satisfies the (Fibered) Isomorphism Conjecture
(in the range $\le N$)
with respect to the family $\calf'(G)$
if and only if $G$ satisfies the (Fibered) Isomorphism Conjecture
(in the range $\le N$) with respect to the family $\calvcyc(G) = \{H \subseteq G\mid H \in
\calvcyc\}$.
\end{theorem}
\begin{proof}
Consider $V \in\calvcyc$. Because of Theorem~\ref{the: transitivity} we have to show that
$V$ satisfies  the Fibered Isomorphism Conjecture
(in the range $\le N$) for the family
$\calf'(V) = \calf'(G) \cap V$. If $V$ is finite, then the claim follows from the
assumptions. It remains to treat the case, where $V$ can be written as an extension
$1 \to \bbZ \to V \xrightarrow{p} F \to 1$ for a finite group $F$.

Since $F$ satisfies the Fibered Isomorphism Conjecture (in the range $\le N$) for
$\calf$ by assumption, $V$ satisfies the Fibered Isomorphism Conjecture (in the range $\le N$) for
$p^*\calf$ by Lemma~\ref{lem: basic inheritance property of fibered conjecture}.
Obviously $p^*\calf(F) \subseteq \calf'(V)$ but in general these two families of
subgroups of $V$ do not agree. Because of
Lemma~\ref{lem: enlarging the family for the Fibered Isomorphism Conjecture}
$V$ satisfies  the Fibered Isomorphism Conjecture (in the range $\le N$) for
$\calf'(V)$.
\end{proof}

The Fibered Isomorphism Conjecture is also well behaved with respect to finite intersections
of families of subgroups.

\begin{lemma} \label{lem: intersections of families}
Let $G$ be a groups and $\calf$ and $\calg$ families of subgroups. Suppose that
$G$ satisfies the Fibered Isomorphism Conjecture (in the range $n \le N$) for both $\calf$ and $\calg$.

Then $G$ satisfies the Fibered Isomorphism Conjecture (in the range $n \le N$) for $\calf \cup \calg$ and
the Fibered Isomorphism Conjecture (in the range $n \le N-1$) for $\calf \cap \calg$.
\end{lemma}

\begin{proof}
We conclude  from Lemma~\ref{lem: enlarging the family for the Fibered Isomorphism Conjecture}
that $G$ satisfies the Fibered Isomorphism Conjecture (in the range $n \le N$) for $\calf \cup \calg$.
Next we prove the claim for $\calf \cap \calg$.

Consider a group homomorphism $\phi \colon K \to G$.
Choose $G$-$CW$-models $\EFG{\calf \cap \calg}{G}$, $\EFG{\calf}{G}$ and
$\EFG{\calg}{G}$ such that $\EFG{\calf \cap \calg}{G}$ is a $G$-$CW$-subcomplex of
both  $\EFG{\calf}{G}$ and $\EFG{\calg}{G}$. This can be arranged by mapping cylinder constructions.
Define a $G$-$CW$-complex
$$X = \EFG{\calf}{G} \cup_{\EFG{\calf \cap \calg}{G}} \EFG{\calg}{G}.$$
For any subgroup $H \subseteq G$ we get
$$X^H = (\EFG{\calf}{G})^H \cup_{(\EFG{\calf \cap \calg}{G})^H} (\EFG{\calg}{G})^H.$$
If $(\EFG{\calf}{G})^H$ and $(\EFG{\calg}{G})^H$ are empty,
the same is true for $X^H$.
If $(\EFG{\calf}{G})^H$ is empty, then
$(\EFG{\calg}{G})^H = X^H$.
If $(\EFG{\calg}{G})^H$ is empty, then
$(\EFG{\calf}{G})^H = X^H$.
If $(\EFG{\calf}{G})^H$, $(\EFG{\calg}{G})^H$ and $(\EFG{\calf \cap \calg}{G})^H$ are contractible,
the same is true for $X^H$. Hence $X$ is a model for $\EFG{\calf \cup \calg}{G}$.
If we apply restriction with $\phi$, we get a decomposition of
$\EFG{\phi^*\calf \cup \phi^*\calg}{K} = \phi^*\EFG{\calf \cup \calg}{G}$
as the union of $\EFG{\phi^*\calf}{K} = \phi^*\EFG{\calf}{G}$ and
$\EFG{\phi^*\calg}{K} = \phi^*\EFG{\calg}{G}$ such that the intersection of
$\EFG{\phi^*\calf}{K}$ and $\EFG{\phi^*\calg}{K}$ is
$\EFG{\phi^*(\calf \cap \calg)}{K} = \phi^*\EFG{\calf \cap \calg}{G}$.
Now the claim follows from the Mayer-Vietoris sequence for $\EFG{\phi^*\calf \cup \phi^*\calg}{K}$
and the Five-Lemma.
\end{proof}

\typeout{--------------------------  Section 2  ------------------------}

\section{General Induction Theorems}
\label{sec: General Induction Theorems}

Let $G$ be a finite group. Let $G\text{-}\FSETS$ be the category of finite $G$-sets.
Morphisms are $G$-maps. Let $\Lambda$ be an associative  commutative ring with unit.
Denote by $\Lambda\text{-}\MOD$ the abelian category of $\Lambda$-modules.
A \emph{bi-functor} $M$  from $G\text{-}\FSETS$ to $\Lambda\text{-}\MOD$
consists of a covariant  functor
$$M_* \colon G\text{-}\FSETS \to \Lambda\text{-}\MOD$$
and a contravariant functor
$$M^* \colon G\text{-}\FSETS \to \Lambda\text{-}\MOD$$
which agree on objects.
\begin{definition}[Mackey functor]
\label{def: Mackey functor}
A \emph{Mackey functor $M$ for $G$ with values in $\Lambda$-modules} is a bifunctor
from $G\text{-}\FSETS$ to $\Lambda\text{-}\MOD$ such that
\begin{itemize}
\item Double Coset formula\\[1mm]
For any cartesian square of finite $G$-sets
\comsquare{S}{\overline{f_1}}{S_1}{\overline{f_2}}{f_1}{S_2}{f_2}{S_0}
the diagram
$$\begin{CD}
M(S) @ > M_*(\overline{f_1}) >> M(S_1)
\\
@A M^*(\overline{f_2}) AA @A  M^*(f_1) AA
\\
M(S_2) @ > M_*(f_2) >> M(S_0)
\end{CD}$$
commutes.

\item Additivity\\[1mm]
Consider two finite $G$-sets $S_0$ and $S_1$. Let
$i_k \colon S_k \to S_0 \amalg S_1$ for $k = 0,1$ be the inclusion. Then the map
$$M^*(i_0) \times M^*(i_1) \colon M(S_0 \amalg S_1) \to  M(S_0) \times M(S_1)$$
is bijective.

\end{itemize}
\end{definition}
One easily checks that the condition Additivity is equivalent to the
requirement that
$$M_*(i_0) \oplus M_*(i_1) \colon M(S_0) \oplus M(S_1) \to M(S_0 \amalg S_1) $$
is bijective since the double coset formula implies
that $(M^*(i_0) \times  M^*(i_1)) \circ (M_*(i_0) \oplus  M_*(i_1))$ is
the identity.

Let $M$ be a Mackey functor. Let $S$ be a finite non-empty $G$-set. Then we get another Mackey functor
$M_S$ by $(M_S)_*(T) = M_*(S \times T)$ and $(M_S)^*(T) = M^*(S \times T)$. The projection
$\pr \colon S \times T \to T$ defines natural transformations of bifunctors from
from $G\text{-}\FSETS$ to $\Lambda\text{-}\MOD$
\begin{eqnarray*}
\theta_S \colon M_S \to M, & & \theta_S(T) = M_*(\pr);
\\
\theta^S \colon M \to M_S, & & \theta^S(T) = M^*(\pr).
\end{eqnarray*}
We call $M$ \emph{$S$-projective} if $\theta_S$ is split surjective as
a natural transformation of bifunctors.
We call $M$ \emph{$S$-injective} if $\theta^S$ is split injective as a natural transformation of bifunctors.

Given a Mackey functor $M$, we can associate to it the $\Lambda$-chain complex
\begin{eqnarray}
\ldots \xrightarrow{c_4} M(S^3) \xrightarrow{c_3}
M(S^2) \xrightarrow{c_2} M(S^1) \xrightarrow{c_1} M(S^0) \to 0 \to 0 \to \ldots
&&
\label{M_*(S^*)}
\end{eqnarray}
and the $\Lambda$-cochain complex
\begin{eqnarray}
\ldots \xleftarrow{c^3} M(S^3) \xleftarrow{c^2}
M(S^2) \xleftarrow{c^1} M(S^1) \xleftarrow{c^0} M(S^0) \leftarrow 0 \leftarrow 0 \leftarrow \ldots.
&&
\label{M^*(S^*)}
\end{eqnarray}
Here $S^n = \prod_{i=1}^n S$ for $n \ge 1$ and $S^0 = G/G$. If $\pr_i^n \colon S^n \to S^{n-1}$
for $i = 1,2, \ldots n$  is the projection omitting the $i$-th coordinate, then
$c_n = \sum_{i=1}^{n} (-1)^i \cdot M_*(\pr_i^n)$ and
$c^n = \sum_{i=1}^{n+1} (-1)^i \cdot M^*(\pr_i^{n+1})$.

The elementary proof of the following result can be found for instance in
\cite[Proposition~6.1.3 and Proposition~6.1.6]{Dieck(1979)}.
The original sources are the papers of Dress \cite {Dress(1973)} and \cite{Dress(1975)}.

\begin{lemma} \label{lem: injectivity and projectivity of Mackey functors}
Let $M$ be a Mackey functor for the finite group $G$ with values in $\Lambda$-modules
and let $S$ be a finite $G$-set.  Then

\begin{enumerate}

\item \label{lem: injectivity and projectivity of Mackey functors: inject. and project.}

The following statements are equivalent:

\begin{enumerate}

\item[(i)] $M$ is $S$-projective;

\item[(ii)] $M$ is $S$-injective;

\item[(iii)] $M$ is a direct summand of $M_S$ as a bifunctor
from $G\text{-}\FSETS$ to $\Lambda\text{-}\MOD$;

\end{enumerate}

\item \label{lem: injectivity and projectivity of Mackey functors: M_S always inj. and proj.}

$M_S$ is always $S$-injective and $S$-projective;

\item \label{lem: injectivity and projectivity of Mackey functors: (co)homology of M(S^*)}

If $M$ is $S$-projective, then the $\Lambda$-chain complex \eqref{M_*(S^*)} is contractible.

If $M$ is $S$-injective, then the $\Lambda$-cochain complex \eqref{M^*(S^*)} is contractible.

\end{enumerate}
\end{lemma}

Let $\Lambda$ be an associative commutative
ring with unit. Let $\Lambda\text{-}\MOD$ be the abelian category of $\Lambda$-modules.
For a small category $\calc$ let $\Lambda\calc\text{-}\MOD$
be the abelian category of contravariant $\Lambda\calc$-modules,
i.e.\ objects are contravariant functors $\calc \to \Lambda\text{-}\MOD$ and morphisms are natural transformations.
In the sequel we use the notation of \cite[Sections~9 and 17]{Lueck(1989)},
where more information about the abelian category
$\Lambda\calc\text{-}\MOD$ and its homological algebra  can be found.
If $T$ is a set, then $\Lambda T$ or $\Lambda (T)$ denote the free $\Lambda$-module with the set $T$ as basis.
The \emph{orbit category} $\Or(G)$ is the small category whose objects are homogeneous $G$-spaces
$G/H$ and whose morphisms are $G$-maps.

\begin{lemma} \label{lem: Tor and Ext of a Mackey functor}
Let $M$ be a Mackey functor for the finite group $G$ with values in $\Lambda$-modules
and let $S$ be a finite non-empty $G$-set. Suppose that $M$ is $S$-injective, or, equivalently,
that $M$ is $S$-projective.  Let $\calf(S)$  be the family of subgroups of $H \subseteq G$
with $S^H \not= \emptyset$. Let $\underline{\Lambda}_{\calf(S)}$ be the contravariant $\Lambda\Or(G)$-module
which sends $G/H$ to $\Lambda$ if $H$ belongs to $\calf(S)$ and to zero otherwise and which sends
a morphism $G/H \to G/K$ to $\id\colon \Lambda \to \Lambda$ if both $H$ and $K$ belong to $\calf(S)$ and to the zero map otherwise.

Then there are natural $\Lambda$-isomorphisms
\begin{eqnarray*}
\Tor_p^{\Lambda\OrF{\calf}{G}}(\underline{\Lambda}_{\calf(S)},M) & \xrightarrow{\cong} &
\left\{\begin{array}{lll}
M(G/G) & & p = 0;
\\
0      & & p \ge 1.
\end{array}\right.
\\
\Ext^p_{\Lambda\OrF{\calf}{G}}(\underline{\Lambda}_{\calf(S)},M) & \xrightarrow{\cong} &
\left\{\begin{array}{lll}
M(G/G) & & p = 0;
\\
0      & & p \ge 1.
\end{array}\right.
\end{eqnarray*}
\end{lemma}
\begin{proof}
Let $T$ be a set, possibly empty. We can assign to it a $\Lambda$-chain
complex $C_*(T)$ which has as $n$-th $\Lambda$-chain module
$\Lambda(T^{n+1})$ and whose $n$-th differential $c_n \colon \Lambda(T^{n+1}) \to \Lambda(T^n)$ is
$\sum_{i=1}^{n+1} (-1)^i \cdot \Lambda(\pr^{n+1}_i)$ for
$\pr_i^{n+1} \colon T^{n+1} \to T^n$ the projection given by leaving out the $i$-th coordinate.
There are in $T$ natural isomorphisms
$$
H_p(C_*(T)) ~ \xrightarrow{\cong} ~
\left\{
\begin{array}{lll}
\Lambda & & p = 0 \text{ and } T \not= \emptyset;
\\

0 & & p \ge 1 \text{ or } T = \emptyset.
\\
\end{array}
\right.
$$
Of course the isomorphism $H_0(C_*(T)) \xrightarrow{\cong} \Lambda$ is induced by the augmentation
$\Lambda$-homomorphism $\Lambda[T] \to \Lambda$ for $T \not= \emptyset$. Fix an element
$t_0 \in T$. The maps
$$\Lambda(T^n) \to \Lambda(T^{n+1}),\quad  (t_1,t_2, \ldots t_n)
~\mapsto ~ -(t_0,t_1,t_2, \ldots t_n)$$
yields the necessary chain homotopies.

Actually, this chain complex
belongs to the fat realization of the nerve of the category which has
$T$ as set of objects and for which the set of morphisms between two
objects consists of exactly one element. Notice that every object
$t_0$ is both initial and terminal.

Next we define a $\Lambda\Or(G)$-chain complex $P_*$ by assigning to $G/H$ the $\Lambda$-chain complex
$C_*(\map_G(G/H,S))$ for the given finite $G$-set $S$. We obtain  a natural identification
of $\Lambda\Or(G)$-modules
$$
H_0(C_*) ~ \xrightarrow{\cong} ~ \underline{\Lambda}_{\calf(S)}.
$$
Obviously each $\Lambda\Or(G)$-chain module $P_n$ is finitely generated
free and hence in particular projective.
Hence $P_*$ is a projective $\Lambda\Or(G)$-resolution of $\underline{\Lambda}_{\calf(S)}$. Thus we get by definition
\begin{eqnarray*}
\Tor_p^{\Lambda\OrF{\calf}{G}}(\underline{\Lambda}_{\calf(S)},M) & = &
H_p\left(P_* \otimes_{\Lambda\Or(G)}M_*\right);
\\
\Ext^p_{\Lambda\OrF{\calf}{G}}(\underline{\Lambda}_{\calf(S)},M) & =&
H^p\left(\Hom_{\Lambda\Or(G)}(P_*,M^*)\right).
\end{eqnarray*}
There is an obvious identification of $\Lambda$-chain complexes
$P_* \otimes_{\Lambda\Or(G)} M_*$ with the $\Lambda$-chain complex
\eqref{M_*(S^*)} if we replace $M(S^0)$  by zero and put $M(S^n)$ as
$(n-1)$-th $\Lambda$-chain module for $n \ge 1$. Analogously
there is an obvious identification of $\Lambda$-chain complexes
$\Hom_{\Lambda\Or(G)}(P_*,M^*)$ with the $\Lambda$-cochain complex
\eqref{M^*(S^*)} if we replace $M(S^0)$  by zero and put $M(S^n)$ as
$(n-1)$-th $\Lambda$-cochain module for $n \ge 1$. Now
Lemma~\ref{lem: Tor and Ext of a Mackey functor} follows from
Lemma~\ref{lem: injectivity and projectivity of Mackey functors}
\eqref{lem: injectivity and projectivity of Mackey functors: (co)homology of M(S^*)}
\end{proof}

Next we deal with the question whether a given Mackey functor is
$S$-projective or, equivalently, $S$-injective.

Let $M$, $N$ and $L$ be bi-functors for the finite group $G$ with
values in $\Lambda$-modules. A pairing
$$M \times N \to L$$
is a family of $\Lambda$-bilinear maps
$$\mu(S) \colon M(S) \times N(S) \to L(S), \quad (m,n) ~ \mapsto ~ \mu(m,n) = m \cdot n$$
indexed by the objects $S$ of $G\text{-}\FSETS$ such that for every morphism $f \colon S \to
T$ in $G\text{-}\FSETS$ we have
\begin{eqnarray}
&\begin{array}{lllll}
L^*(f)(x \cdot y) & = & M^*(f)(x) \cdot N^*(f)(y), & & x \in M(T), y \in N(T);
\\
x \cdot N_*(f)(y) & = & L_*(f)(M^*(f)(x) \cdot y),  & & x \in M(T), y \in N(S);
\\
M_*(f)(x) \cdot y & = & L_*(f)(x \cdot N^*(f)(y)),  & & x \in M(S), y \in N(T).
\end{array} &
\label{conditions for a pairing}
\end{eqnarray}

\begin{definition}[Green functor] A \emph{Green functor} for the finite group $G$ with values in $\Lambda$-modules
is a Mackey functor $U$ together with a pairing
$$\mu \colon U \times U \to U$$
and a choice of elements $1_S \in U(S)$ for each finite $G$-set $S$
such that for each finite $G$-set $S$ the pairing
$\mu(S) \colon U(S) \times U(S) \to U(S)$ and the element $1_S$ determine the structure of
an associative $\Lambda$-algebra with unit on $U(S)$.
Moreover, it is required that $U^*(f)(1_T) = 1_S$ for every morphism $f \colon S \to
T$ in $G\text{-}\FSETS$.

A (left) $U$-module $M$ is a Mackey functor for the finite group $G$ with values in
$\Lambda$-modules together with a pairing
$$\nu \colon U \times M \to M$$
such that for every finite $G$-set $S$ the pairing $\nu(S) \colon U(S) \times M(S) \to
M(S)$ defines the structure of a $U(S)$-module on $M(S)$,
where $1_S$ acts as $\id_{M(S)}$.
\end{definition}

The proof of the next result can be found for instance in
\cite[Theorem~6.2.2]{Dieck(1979)}.

\begin{theorem}[Criterion for $S$-projectivity]
\label{the: Criterion for S-projectivity}
Let $G$ be a finite group and let $\Lambda$ be an associative commutative ring with unit. Let $S$ be a finite
$G$-set. Let $U$ be a Green
functor. Then the following statements are equivalent:

\begin{enumerate}

\item The projection $\pr \colon S \to G/G$ induces an epimorphism
$U_*(\pr)\colon U(S) \to U(G/G)$;

\item $U$ is $S$-injective;

\item Every $U$-module is $S$-injective.

\end{enumerate}

\end{theorem}

The results of this section have the following application to equivariant homology theories.

\begin{theorem}[Criterion for induction for $G$-homology theories]
\label{the: Criterion for induction for equivariant homology theories}
Let $G$ be a finite group and $\calf$
be a family of subgroups of $G$.
Let $\calh^G_*$ be a $G$-homology theory with values in $\Lambda$-modules for
an associative commutative ring $\Lambda$ with unit.
Suppose that the following conditions are satisfied:

\begin{itemize}

\item There exists a Green functor $U$ for the finite group $G$ with values in $\Lambda$-modules
such that the $\Lambda$-homomorphism
$$\bigoplus_{H \in \calf} U(\pr_H) \colon \bigoplus_{H \in \calf} U(G/H) \to U(G/G)$$
is surjective, where $\pr_H \colon G/H \to G/G$ is the projection.

\item For every $n \in \bbZ$ there is a (left) $U$-module $M$ such that the covariant functor
$M_* \colon G\text{-}\FSETS \to \Lambda\text{-}\MOD$ is naturally equivalent to the
covariant functor
$$
\calh^G_n\colon G\text{-}\FSETS \to \Lambda\text{-}\MOD, \quad S ~ \mapsto ~ \calh^G_n(S).
$$

\end{itemize}

Then the projection $\pr \colon \EFG{\calf}{G} \to G/G$ induces for all $n \in \bbZ$ a
$\Lambda$-isomorphism
\begin{equation} \label{eq: homology iso}
\calh^G_n(\pr) \colon \calh^G_n(\EFG{\calf}{G}) \xrightarrow{\cong} \calh_n^G(G/G).
\end{equation}
and the canonical map
\begin{equation} \label{eq: colim iso}
\colim_{\OrF{\calf}{G}} \calh^G_n(G/?) ~ \xrightarrow{\cong} \calh^G_n(G/G)
\end{equation}
is bijective, where $\OrF{\calf}{G} \subset \Or(G)$ is the full subcategory of the orbit category
whose objects are homogeneous spaces $G/H$ with $H \in \calf$.
\end{theorem}

\begin{proof}
Let $S$ be the $G$-set $\amalg_{H \in \calf} G/H$.
The first condition together with
Lemma \ref{lem: injectivity and projectivity of Mackey functors}
\eqref{lem: injectivity and projectivity of Mackey functors:  inject. and project.}
and Theorem~\ref{the: Criterion for S-projectivity} implies
that $M$ is $S$-projective. We conclude from the second condition and
Lemma~\ref{lem: Tor and Ext of a Mackey functor} that there is a canonical $\Lambda$-isomorphism
\begin{eqnarray*}
\Tor_p^{\Lambda\OrF{\calf}{G}}(\underline{\Lambda}_{\calf(S)},\calh^G_q(G/?)) & \xrightarrow{\cong} &
\left\{\begin{array}{lll}
\calh^G_q(G/G) & & p = 0;
\\
0      & & p \ge 1,
\end{array}\right.
\end{eqnarray*}
for all $q \in \bbZ$. The cellular $\Lambda\Or(G)$-chain complex of $\EFG{\calf}{G}$ is a
projective $\Lambda\Or(G)$-resolution of $\underline{\Lambda}_{\calf(S)}$.
Hence $\Tor_p^{\Lambda\OrF{\calf}{G}}(\underline{\Lambda}_{\calf(S)},\calh^G_q(G/?)$ agrees
with the Bredon homology $H_p^{\Lambda\Or(G)}(\EFG{\calf}{G};\calh^G_q(G/?))$. But this is
exactly the $E^2$-term in the equivariant Atiyah-Hirzebruch spectral sequence
which converges to $\calh^G_{p+q}(\EFG{\calf}{G})$. This implies that the $E^2$-term is concentrated
on the $y$-axis and the spectral sequence collapses. Hence the edge homomorphism yields an isomorphism
$$\calh^G_n(\EFG{\calf}{G}) \xrightarrow{\cong} \calh_n^G(G/G)$$
for all $n \in \bbZ$ which can easily be identified with the $\Lambda$-map $\calh^G_n(\pr)$.
There is a natural identification
\begin{eqnarray*}
\Tor_0^{\Lambda\OrF{\calf}{G}}(\underline{\Lambda}_{\calf(S)},\calh^G_q(G/?)) & \xrightarrow{\cong} &
\colim_{\OrF{\calf}{G}} \calh^G_q(G/?).
\end{eqnarray*}
This finishes the proof of Theorem~\ref{the: Criterion for induction for equivariant homology theories}.
\end{proof}

\begin{remark} \em
The bijectivity of the map \eqref{eq: colim iso} is a consequence of the
exactness of the complex \eqref{M_*(S^*)} at $M(S^0)$ and $M(S^1)$.
It is not hard to see that this map factors as
\begin{equation*}
\colim_{\OrF{\calf}{G}} \calh^G_n(G/?) \to \calh^G_n(\EFG{\calf}{G}) \to \calh_n^G(G/G)
\end{equation*}
and therefore \eqref{eq: homology iso} is onto.
In order to get our applications to assembly maps it is important that \eqref{eq: homology iso}
is in fact an isomorphism.
(Surjectivity alone is not helpful, since the Transitivity Principle~\ref{the: transitivity}
does not apply to surjections.) The proof of the bijectivity of \eqref{eq: homology iso}
is based on the result due to Dress that the complex \eqref{M_*(S^*)} is not only exact at $M(S^0)$ and $M(S^1)$ but is contractible.

The dual version of the map \eqref{eq: colim iso} leads to the following induction result in
\cite{Farrell-Hsiang(1977)} that for a group $\Gamma$ (in their setting a Bieberbach group)
 which can be written as an extension
$1 \to \Gamma_0 \to \Gamma \xrightarrow{p} G \to 1$ for a finite group $G$ the map
$$L_n(\bbZ\Gamma) \otimes_{\bbZ} \bbQ \xrightarrow{\cong} \lim L_n(\bbZ[p^{-1}(C)]) \otimes_{\bbZ} \bbQ$$
is bijective, where the inverse limit runs over the cyclic subgroups of $G$. We expect that there is a dual
version of \eqref{eq: homology iso} which leads to kind of dual version of induction theorems for groups mapping surjectively
to finite groups.
\em
\end{remark}

\typeout{--------------------------  Section 3  ------------------------}

\section{The Swan Group as a Functor on Groupoids}
\label{sec: The Swan Group as a Functor on Groupoids}

In this section we will construct Green functors and modules over them that will be used
in the application of Theorem~\ref{the: Criterion for induction for equivariant homology theories}
in the proof of Theorem~\ref{the: Induction theorem for algebraic K-theory} in
Section~\ref{sec: Proof of the Main Result for Algebraic $K$-Theory}.

Let $\calg$ be a small groupoid.
Our main example is the {\em transport groupoid} $\calg^G(S)$ of a $G$-set $S$.
The set of objects is given by the set $S$ itself.
The set of morphisms from $s_1$ to $s_2$ consists of those elements $g \in G$ which satisfy
$gs_1 = s_2$. Composition comes from the multiplication in $G$. A $G$-map
$f \colon S \to T$ induces a functor $\calg^G(f) \colon \calg^G(S) \to \calg^G(T)$.
Let $R$ be an associative ring with unit.
Denote by
$R\calg\text{-}\FGPMOD$ the category of contravariant functors from $\calg$ to finitely
generated projective $R$-modules and by $K_n(R\calg)$ its $K$-theory.
Let $\func(\calg,\bbZ\text{-}\FGMOD)$ and
$\func(\calg,\bbZ\text{-}\FGFMOD)$ respectively be the category of contravariant functors from $\calg$
to the category $\bbZ\text{-}\FGMOD$ of finitely generated $\bbZ$-modules and
to the category $\bbZ\text{-}\FGFMOD$ of finitely generated free $\bbZ$-modules
respectively. Let $\Sw(\calg)$ and $\Sw^f(\calg)$ respectively be the $K_0$-group
of $\func(\calg,\bbZ\text{-}\FGMOD)$ and
$\func(\calg,\bbZ\text{-}\FGFMOD)$ respectively. The forgetful map
$$\Sw^f(\calg) \xrightarrow{\cong} \Sw(\calg)$$
is a bijection. This is proved for groups in \cite[page 890]{Pedersen-Taylor(1978)}
and carries easily over to groupoids.
We will concentrate our discussion mostly on $\Sw^f(\calg)$.

Given a contravariant $\bbZ\calg$-module
$M$ and a contravariant $R\calg$-module $N$, let $M \otimes_{\bbZ} N$ be the  contravariant
$R\calg$-module which assigns to an object $c$ the $R$-module $M(c) \otimes_{\bbZ} N(c)$.
If $M$ belongs to $\func(\calg,\bbZ\text{-}\FGFMOD)$, then the functor
$M \otimes_{\bbZ} -$ is exact. If $M$ and $N$ belong to $\func(\calc,\bbZ\text{-}\FGFMOD)$, then
$M \otimes_{\bbZ} N$ belongs to $\func(\calc,\bbZ\text{-}\FGFMOD)$.
Hence $\otimes_{\bbZ}$ induce a pairing
\begin{eqnarray}
\label{eq: ring-structure-on-Sw^f}
\Sw^f(\calg) \times \Sw^f(\calg) \to \Sw^f(\calg).
\end{eqnarray}
With this pairing $\Sw^f(\calg)$ is a commutative associative ring
with the class of the constant contravariant $\bbZ\calc$-module with value $\bbZ$ as unit.

If $M$ belongs to $\func(\calg,\bbZ\text{-}\FGFMOD)$,
then $M \otimes_{\bbZ} -$ sends finitely generated projective
$R\calg$-modules to finitely generated projective $R\calg$-modules. This is well-known if
the groupoid $\calc$ has only one object, i.e.\ in case of group rings, and hence holds
also for a groupoid $\calg$. Hence we get for a groupoid a pairing
\begin{eqnarray}
\label{eq: module-structure-on-K}
\Sw^f(\calg) \times K_n(R\calg) \to K_n(R\calg)
\end{eqnarray}
which turns $K_n(R\calg)$ into a $\Sw^f(\calg)$-module.

Given a functor $F \colon \calc \to \cald$, induction defines an functor
$$RF_* \colon R\calc\text{-}\MOD \to R\cald\text{-}\MOD$$
and restriction defines a functor
$$RF^* \colon R\cald\text{-}\MOD\to  R\calc\text{-}\MOD.$$
The induction functor is given by
$- \otimes_{R\calc} R\mor_{\cald}(??,F(?))$ and restriction by
$- \otimes_{R\cald} R\mor_{\cald}(F(?),??)$ for the
$R\calc\text{-}R\cald$-bimodule  $R\mor_{\cald}(??,F(?))$ and the
$R\cald\text{-}R\calc$-bimodule  $ R\mor_{\cald}(F(?),??)$,
where $?$ runs through the objects in $\calc$ and $??$ through the objects in $\cald$
 (see \cite[9.15 on page 166]{Lueck(1989)}).
We have the obvious equalities $(F_2 \circ F_1)_* = (F_2)_* \circ (F_1)_*$ and
$(F_2 \circ F_1)^* = (F_1)^* \circ (F_2)^*$ for functors
$F_1 \colon \calc_1 \to \calc_2$ and $F_2 \colon \calc_2 \to \calc_3$.

Let $F \colon \calc \to \cald$ be a functor of groupoids. Induction with $F$ sends
finitely generated projective $R\calc$-modules to finitely generated projective
$R\cald$-modules and respects direct sums. Hence induction induces for $n \in \bbZ$
homomorphisms of abelian groups
$$F_* \colon K_n(R\calc) \to K_n(R\cald).$$
Restriction with $F$ yield  exact functors
$\func(\cald,\bbZ\text{-}\FGFMOD) \to \func(\calc,\bbZ\text{-}\FGFMOD)$ and
$\func(\cald,\bbZ\text{-}\FGMOD) \to \func(\calc,\bbZ\text{-}\FGMOD)$
and thus ring homomorphisms
\begin{eqnarray*}
F^* \colon \Sw^f(\cald) & \to & \Sw^f(\calc);
\\
F^* \colon \Sw(\cald) & \to & \Sw(\calc).
\end{eqnarray*}
We call $F$ \emph{admissible} if for each object $c \in \calc$ the group
homomorphism $\aut_{\calc}(c) \to \aut_{\cald}(F(c))$ induced by $F$ is injective
and its image has finite index and the map $\pi_0(F) \colon \pi_0(\calc) \to \pi_0(\cald)$
has the property that the preimage of any element in $\pi_0(\cald)$ is finite.
Note that if $G$ is a (not necessary finite) group and
$f : S \to T$ is a map of finite $G$-sets then $\calg^G(f)$ is admissible.
For admissible $F$ induction and restriction do also induce
homomorphisms of abelian groups
\begin{eqnarray*}
F^* \colon K_n(R\cald) & \to & K_n(R\calc);
\\
F_* \colon \Sw^f(\calc) & \to & \Sw^f(\cald);
\\
F_* \colon \Sw(\calc) & \to & \Sw(\cald).
\end{eqnarray*}
The various claims above are well-known for groups, i.e.\ groupoids with one object
and therefore carry easily over to groupoids.

Let $E,F \colon \calc \to \cald$ be functors which are naturally equivalent. Then
we get the following equalities of homomorphisms:
\begin{eqnarray*}
E^* = F^* \colon \Sw^f(\cald) & \to & \Sw^f(\calc);
\\
E^* = F^* \colon \Sw(\cald) & \to & \Sw(\calc);
\\
E_* = F_* \colon K_n(R\calc)      & \to & K_n(R\cald)
\end{eqnarray*}
and $E$ is admissible if and only if $F$ is and in this case also the following
homomorphisms agree
\begin{eqnarray*}
E_* = F_* \colon \Sw^f(\calc) & \to & \Sw^f(\cald);
\\
E^* = F^* \colon \Sw(\calc) & \to & \Sw(\cald);
\\
E^* = F^* \colon K_n(R\cald)      & \to & K_n(R\calc).
\end{eqnarray*}

Let $\calg$ be a groupoid with a finite set $\pi_0(\calg)$ of components.
For a component $\calc \in \pi_0(\calg)$ let $i_{\calc} \colon \calc \to \calg$ be the inclusion.
Then we obtain to one another inverse isomorphisms
\begin{eqnarray*}
\bigoplus_{\calc \in \pi_0(\calg)} (i_{\calc})_* \colon
\bigoplus_{\calc \in \pi_0(\calg)} K_n(R\calc) & \xrightarrow{\cong} &
K_n(R\calg);
\\
\prod_{\calc \in \pi_0(\calg)} (i_{\calc})^* \colon  K_n(R\calg) & \xrightarrow{\cong}  &
\prod_{\calc \in \pi_0(\calg)} K_n(R\calc)
\end{eqnarray*}
and to one another inverse isomorphisms
\begin{eqnarray*}
\bigoplus_{\calc \in \pi_0(\calg)} (i_{\calc})_* \colon
\bigoplus_{\calc \in \pi_0(\calg)} \Sw^f(\calc) & \xrightarrow{\cong} &
\Sw^f(\calg);
\\
\prod_{\calc \in \pi_0(\calg)} (i_{\calc})^* \colon  \Sw^f(\calg) & \xrightarrow{\cong}  &
\prod_{\calc \in \pi_0(\calg)} \Sw^f(\calc)
\end{eqnarray*}
and similarly for $\Sw$. If $\calg$ is a connected groupoid and $x$ an object in $\calg$,
then the inclusion $i \colon \aut_{\calg}(x) \to \calg$ induces to another inverse
isomorphisms
\begin{eqnarray*}
i_* \colon K_n(R[\aut_{\calg}(x)]) & \xrightarrow{\cong} & K_n(R\calg);
\\
i^* \colon K_n(R\calg) & \xrightarrow{\cong} & K_n(R[\aut_{\calg}(x)]),
\end{eqnarray*}
and analogously for $\Sw^f$ and $\Sw$.

\begin{lemma} \label{lem_ Mackey for calg^g applied to a pullback}
Consider a cartesian square of $G$-sets
\comsquare{S}{\overline{f_1}}{S_1}{\overline{f_2}}{f_1}{S_2}{f_2}{S_0}
Then the following diagram of functors commutes up to natural
equivalence
$$\begin{CD}
R\calg^G(S)\text{-}\MOD @ > R\calg^g(\overline{f_1})_* >> R\calg^G(S_1)\text{-}\MOD
\\
@A  R\calg^G(\overline{f_2})^* AA @A  R\calg^G(f_1)^* AA
\\
R\calg^G(S_2)\text{-}\MOD @ > R\calg^g(f_2)_* >> R\calg^G(S_0)\text{-}\MOD
\end{CD}$$
\end{lemma}
\begin{proof}
The composition
$$R\calg^G(S_2)\text{-}\MOD \xrightarrow{R\calg^g(f_2)_*} R\calg^G(S_0)\text{-}\MOD
\xrightarrow{R\calg^G(f_1)^*} R\calg^G(S_1)\text{-}\MOD$$
is given by
$$- \otimes_{R\calg^G(S_2)} R\mor_{\calg^G(S_0)}(!,\calg^G(f_2)(??))
\otimes_{R\calg^G(S_0)} R\mor_{\calg^G(S_0)}(\calg^G(f_1)(?),!),$$
 where $?$, $??$ and $!$ run through objects in
$\calg^G(S_1)$, $\calg^G(S_2)$ and $\calg^G(S_0)$.
The composition
$$R\calg^G(S_2)\text{-}\MOD \xrightarrow{R\calg^g(\overline{f_2})^*} R\calg^G(S)\text{-}\MOD
\xrightarrow{R\calg^G(\overline{f_1})_*} R\calg^G(S_1)\text{-}\MOD$$
is given by
$$- \otimes_{R\calg^G(S_2)} R\mor_{\calg^G(S_2)}(\calg^G(\overline{f_2})(!!),??)
\otimes_{R\calg^G(S)}R\mor_{\calg^G(S_1)}(?,\calg^G(\overline{f_1})(!!)),$$
 where $?$, $??$ and $!!$ run through objects in
$\calg^G(S_1)$, $\calg^G(S_2)$ and $\calg^G(S)$. Hence it suffices to construct an isomorphism
of $R\calg^G(S_2)$-$R\calg^G(S_1)$-bimodules
\begin{multline*}
R\mor_{\calg^G(S_0)}(!,\calg^G(f_2)(??))
\otimes_{R\calg^G(S_0)} R\mor_{\calg^G(S_0)}(\calg^G(f_1)(?),!)
\\
\xrightarrow{\cong}
R\mor_{\calg^G(S_2)}(\calg^G(\overline{f_2})(!!),??)
\otimes_{R\calg^G(S)}R\mor_{\calg^G(S_1)}(?,\calg^G(\overline{f_1})(!!)).
\end{multline*}
For this purpose it suffices to construct an  isomorphism of
$\calg^G(S_2)$-$\calg^G(S_1)$-bisets
\begin{multline}
\mor_{\calg^G(S_0)}(!,\calg^G(f_2)(??))
\otimes_{\calg^G(S_0)} \mor_{\calg^G(S_0)}(\calg^G(f_1)(?),!)
\\
\xrightarrow{\cong}
\mor_{\calg^G(S_2)}(\calg^G(\overline{f_2})(!!),??)
\otimes_{\calg^G(S)}\mor_{\calg^G(S_1)}(?,\calg^G(\overline{f_1})(!!))
\label{isomorphism of bisets}
\end{multline}
where $\otimes_{\calg^G(S)}$ is now to be understood with respect to
the category of sets.  There is an obvious bijection of
$\calg^G(S_2)$-$\calg^G(S_1)$-bisets
\begin{multline*}
\mor_{\calg^G(S_0)}(!,\calg^G(f_2)(??))
\otimes_{\calg^G(S_0)} \mor_{\calg^G(S_0)}(\calg^G(f_1)(?),!)
\\
\xrightarrow{\cong}
\mor_{\calg^G(S_0)}(\calg^G(f_1)(?),\calg^G(f_2)(??))
\end{multline*}
which sends $u \otimes v$ to $u \circ v$. Its inverse sends
$w \colon \calg^G(f_1)(?) \to \calg^G(\overline{f_2})(??)$
to $\id_{\calg^G(\overline{f_2})(??)} \otimes w$.
There is a map of
$\calg^G(S_2)$-$\calg^G(S_1)$-bisets
\begin{multline*}
\mor_{\calg^G(S_2)}(\calg^G(\overline{f_2})(!!),??)
\otimes_{\calg^G(S)}\mor_{\calg^G(S_1)}(?,\calg^G(\overline{f_1})(!!))
\\
\xrightarrow{\cong}
\mor_{\calg^G(S_0)}(\calg^G(f_1)(?),\calg^G(f_2)(??))
\end{multline*}
which sends $u \otimes v$ to
$\calg^G(f_2)(u) \circ \calg^G(f_1)(v)$. This makes sense, since
$\calg^G(f_1) \circ \calg^G(\overline{f_1})$ and
$\calg^G(f_2) \circ \calg^G(\overline{f_2})$ coincide.
In order to show that this map is a bijection of bisets,
we need the assumption that the square of $G$-sets
appearing in Lemma~\ref{lem_ Mackey for calg^g applied to a pullback} is cartesian.
Namely, we construct the inverse
\begin{multline*}
\mor_{\calg^G(S_0)}(\calg^G(f_1)(?),\calg^G(f_2)(??))
\\
\xrightarrow{\cong}
\mor_{\calg^G(S_2)}(\calg^G(\overline{f_2})(!!),??)
\otimes_{\calg^G(S)}\mor_{\calg^G(S_1)}(?,\calg^G(\overline{f_1})(!!))
\end{multline*}
as follows. Consider a morphism
$w \colon \calg^G(f_1)(?) \to \calg^G(\overline{f_2})(??)$.
It is given by an element $g \in G$ satisfying $g \cdot f_1(?) = f_2(??)$.
The element $(g\cdot ?,??) \in S_1 \times S_2$ defines a unique element $!! \in S$
since $f_1(g \cdot ?) = f_2(??)$. We have $\overline{f_1}(!!) = g \cdot ?$ and
$\overline{f_2}(!!) = ??$. The element $g$ defines a morphism
$u \colon ? \to \calg^G(\overline{f_1})(!!)$. Now define the image of $w$ by
$\id_{??} \otimes u$. One easily checks that these
two maps of bisets are inverse to one another. The desired map of bisets
\eqref{isomorphism of bisets} is given by the two bijections of bisets above. This finishes the proof
of Lemma~\ref{lem_ Mackey for calg^g applied to a pullback}.
\end{proof}

We can now define the Green functors that we will need.
Let $\phi : K \to G$ be a group homomorphism whose target is a finite group $G$.
Then $S \mapsto \Sw^f(\calg^G(S))$ and $S \mapsto \Sw^f(\calg^K(\phi^* S))$ are
Green functors.
(Here $\phi^* S$ denotes the finite $K$-set obtained by restricting
the finite $G$-set with $\phi$.)
The covariant functorial
structure comes from induction and the contravariant functorial structure from
restriction with $\calg^G(f)$ and $\calg^K(\phi^* f)$. The double coset formula
for $\Sw^f(\calg^G(-))$ follows directly from Lemma~\ref{lem_ Mackey for calg^g applied to a pullback}
and for $\Sw^f(\calg^K(\phi^*-)$ from
Lemma~\ref{lem_ Mackey for calg^g applied to a  pullback} using the fact that $\phi^*$
sends a cartesian square of finite $G$-sets to a cartesian square of finite $K$-sets.
The required pairings come from \eqref{eq: ring-structure-on-Sw^f}.
We leave the verification of the
conditions \eqref{conditions for a pairing} to the reader, they follow
from certain natural equivalences of functors on the level of
$R\calg$-modules.

The group homomorphism $\phi$ induces a functor $\calg^{\phi}(S) \colon
\calg^K(\phi^*S) \to \calg^G(S)$ for every finite $F$-set which is natural in
$S$. Restriction with it induces a morphisms of Green functors for $G$
$$\phi^* \colon \Sw^f(\calg^G(-)) \to \Sw^f(\calg^K(\phi^*-).$$

Analogously one constructs a left $\Sw^f(\calg^K(\phi^*-)$-module $K_n(\calg^K(\phi^*-)$.
The Mackey structure on $K_n(\calg^K(\phi^*-)$
comes from induction and restriction and the desired pairing from \eqref{eq: module-structure-on-K}.
Using $\phi^* \colon \Sw^f(\calg^G(-)) \to \Sw^f(\calg^K(\phi^*-)$
we obtain a left $\Sw^f(\calg^G(-))$-module structure on
$K_n(\calg^K(\phi^*-)$.

If $\bbC \subseteq R$, we can replace $\Sw^f(\calg)$ by the version
$\Sw(\bbC \calg)$, where one uses modules over $\bbC$ instead of
$\bbZ$. If $H$ is a finite group, then $\Sw(\bbC H)$ agrees with the
complex representation ring of $H$.

\typeout{--------------------------  Section 4  ------------------------}

\section{Proof of the Main Result for Algebraic $K$-Theory}
\label{sec: Proof of the Main Result for Algebraic $K$-Theory}

Let $R$ be an associative ring with unit.
Let $\calh^?_*(-;\bfK_R)$ be the equivariant homology theory associated to the covariant
functor $\bfK_R$ from the category $\GROUPOIDS$ of small groupoids to the category of
$\Omega\text{-}\SPECTRA$ which sends a groupoid $\calg$ to the (non-connective) algebraic
$K$-theory spectrum of the category of finitely generated projective contravariant  $R\calg$-modules.
Notice that
$$K_n(R\calg) = \pi_n(\bfK_R(\calg))$$
and
$$\calh^G_n(G/H;\bfK_R) ~ = ~ \calh^H_n(\pt;\bfK_R) ~ = ~  K_n(RH)$$
holds for all $n \in \bbZ$. (See \cite[Chapter 6]{Lueck-Reich(2004g)} for more details).
Recall that the Farrell-Jones Conjecture for algebraic
$K$-theory and a group $G$ says that the projection $\pr \colon \EFG{\calvcyc}{G} \to G/G$ induces for all
$n \in \bbZ$ an isomorphism
$$\calh^G_n(\pr;\bfK_R) \colon \calh^G_n(\EFG{\calvcyc}{G};\bfK_R) ~ \xrightarrow{\cong} ~
\calh^G_n(G/G;\bfK_R) = K_n(RG).$$

Let $\Lambda$ be a
commutative associative ring with unit such that $\Lambda$ is flat as a $\bbZ$-module,
or, equivalently, $\Lambda$ is torsionfree as an abelian group.
Let $\phi \colon K \to G$ be a group homomorphism.
Then $X \mapsto \Lambda \otimes_{\bbZ} \calh^K_*(\phi^*X;\bfK_R)$
defines a $G$-homology theory, where for a $G$-$CW$-complex $X$ we denote by $\phi^*X$ the
$K$-$CW$-complex obtained from $X$ by restriction with $\phi$.

\begin{lemma} \label{lem: calh^?_*(-,bfK) and assumptions for induction}
Let $\calf$ be a class of finite groups closed under isomorphism and taking subgroups.
Let $\phi \colon K \to G$ be a group homomorphism with a finite group $G$ as target.

Then the $G$-homology theory
$\Lambda \otimes_{\bbZ} \calh^K_*(\phi^*-;\bfK_R)$
satisfies the assumptions appearing in
Theorem~\ref{the: Criterion for induction for equivariant homology theories} for the
family $\calf(G) = \{H \subseteq G \mid H \in \calf\}$ in the following cases:

\begin{enumerate}

\item $\calf$ is the class $\calh$ of hyperelementary groups and
      $\Lambda = \bbZ$;

\item $\calf$ is the class $\cale$ of elementary groups and
      $\Lambda = \bbZ$, provided $\bbC \subseteq R$;

\item For a given prime $p$ the family $\calf$ is the class $\calh_p$ of
      $p$-hyperelementary groups and $\Lambda = \bbZ_{(p)}$;

\item For a given prime $p$ the family $\calf$ is the class $\cale_p$ of
      $p$-elementary groups and $\Lambda = \bbZ_{(p)}$, provided $\bbC \subseteq R$;

\item $\calf$ is the class $\calfcyc$ of finite cyclic groups and $\Lambda = \bbQ$.

\end{enumerate}

\end{lemma}

\begin{proof}
In Section~\ref{sec: The Swan Group as a Functor on Groupoids} we
have constructed the  Green functor
$\Lambda \otimes \Sw^f(\calg^G(-))$ and the
$\Lambda \otimes \Sw^f(\calg^G(-))$-module
$\Lambda \otimes K_n(R\calg^K(\phi^*-))$.
There is a natural equivalence of covariant functors
$G\text{-}\FSETS \to \bbZ-\MOD$ from
$K_n(R\calg^K(\phi^*-))$ to $\calh^G_n(\phi^*-;\bfK_R)$.
Hence it remains to check the following

\begin{enumerate}

\item \label{swans-induction-for-calh}
      The map coming from
      induction with respect to the various inclusions $H \subseteq G$
      $$\bigoplus_{H \in \calh} \Sw^f(H) ~ \to ~ \Sw^f(G)$$
      is surjective;

\item \label{swans-induction-for-cale}
      The map coming from
      induction with respect to the various inclusions $H \subseteq G$
      $$\bigoplus_{H \in \cale} \Sw(\bbC H) ~ \to ~ \Sw(\bbC G)$$
      is surjective;

\item \label{swans-induction-for-calh_p}
      For a given prime $p$  the map coming from
      induction with respect to the various inclusions $H \subseteq G$
      $$\bigoplus_{H \in \calh_p} \Sw^f(H)_{(p)} ~ \to ~ \Sw^f(G)_{(p)}$$
      is surjective;

\item \label{swans-induction-for-cale_p}
      For a given prime $p$  the map coming from
      induction with respect to the various inclusions $H \subseteq G$
      $$\bigoplus_{H \in \cale_p} \Sw(\bbC H)_{(p)} ~ \to ~ \Sw(\bbC G)_{(p)}$$
      is surjective;

\item \label{swans-induction-for-calfcyc}
      The map coming from
      induction with respect to the various inclusions $H \subseteq G$
      $$\bigoplus_{H \in \calfcyc} \bbQ \otimes_{\bbZ} \Sw^f(H)  ~ \to ~ \bbQ \otimes_{\bbZ}\Sw^f(G)$$
      is surjective.
\end{enumerate}
\eqref{swans-induction-for-calh} and \eqref{swans-induction-for-calfcyc} are
proved by Swan \cite[Corollary 4.2]{Swan(1960a)} for
$\Sw(\calg^G(-))$ which is isomorphic to $\Sw^f(\calg^G(-))$ by
\cite[Proposition 1.1]{Swan(1960a)}.
We conclude from  \cite[Lemma 4.1]{Swan(1960a)} and
\cite[Section 12]{Swan(1963)} that the torsion elements of $\Sw(G)$ are all nilpotent.
Together with \eqref{swans-induction-for-calfcyc} and Theorem~\ref{the: Criterion for S-projectivity}
this show that the assumptions of \cite[6.3.3]{Dieck(1979)} are satisfied.
Thus \eqref{swans-induction-for-calh_p} follows from
\cite[6.3.3]{Dieck(1979)}. We get \eqref{swans-induction-for-cale} and
\eqref{swans-induction-for-cale_p} from
\cite[Theorem 27 and Theorem 28 on page 98]{Serre(1977)}, since $\Sw(\bbC H)$ for a finite group $H$ is the
same as the complex representation ring of $H$.
\end{proof}

Next we give the proof of Theorem~\ref{the: Induction theorem for algebraic K-theory}.

\begin{proof}
If we combine
Theorem~\ref{the: Criterion for induction for equivariant homology theories}
and
Lemma~\ref{lem: calh^?_*(-,bfK) and assumptions for induction},
then we get for any group homomorphism $\phi \colon K \to G$ whose target $G$ is a finite group
that the projection $\pr \colon \phi^*\EFG{\calf(G)}{G} \to \pt$ induces for all $n \in \bbZ$ bijections
\begin{equation*}
\id_{\Lambda} \otimes_{\bbZ} \calh^K_n(\pr;\bfK_R) :
\Lambda \otimes_{\bbZ}  \calh^K_n(\phi^*\EFG{\calf(G)}{G};\bfK_R) ~ \to ~
\Lambda \otimes_{\bbZ} \calh^K_n(\pt;\bfK_R)
\end{equation*}
in the following cases

\begin{enumerate}

\item $\calf$ is the class $\calh$ of hyperelementary groups and
      $\Lambda = \bbZ$;

\item $\calf$ is the class $\cale$ of elementary groups and $\Lambda =
      \bbZ$, provided  that $\bbC \subseteq R$;

\item For a given prime $p$ the family $\calf$ is the class $\calh_p$ of
      $p$-hyperelementary groups and $\Lambda = \bbZ_{(p)}$;

\item For a given prime $p$ the family $\calf$ is the class $\cale_p$ of
      $p$-elementary groups and $\Lambda = \bbZ_{(p)}$, provided that
      $\bbC \subseteq R$;

\item $\calf$ is the class $\calfcyc$ of finite cyclic groups and $\Lambda = \bbQ$.

\end{enumerate}

Since $\phi^*\EFG{\calf(G)}{G} = \EFG{\phi^*\calf(G)}{K}$ holds, this shows
that for the finite group $G$  and the equivariant homology theory
$\Lambda \otimes_{\bbZ} \calh^?_*(-;\bfK_R)$
the Fibered Isomorphism Conjecture
\ref{def: (Fibered) Isomorphism Conjectures for calh^?_*} for the family
$\calf(G)$ holds. Now apply Theorem~\ref{the: calfin to calvcyc} and use the fact that
the (Fibered) Isomorphism Conjecture in the sense of
Definition~\ref{def: (Fibered) Isomorphism Conjectures for calh^?_*}
for $\Lambda \otimes_{\bbZ} \calh^?_*(-;\bfK_R)$ is the same as
the original (Fibered) Farrell-Jones Conjecture for algebraic $K$-theory
for group rings with coefficients in $R$.
For the Farrell-Jones Conjecture this is indicated in \cite[p.239]{Davis-Lueck(1998)},
a careful proof can be found in \cite[Corollary 9.2]{Hambleton-Pedersen(2004)}.
For the Fibered Farrell-Jones Conjecture compare
\cite[Remark 6.6]{Bartels-Lueck(2004)} and
\cite[Remark 4.14]{Lueck-Reich(2004g)}.
This finishes the proof of  assertions
\eqref{the: Induction theorem for algebraic K-theory: calh} and
\eqref{the: Induction theorem for algebraic K-theory: calh_p} of
Theorem~\ref{the: Induction theorem for algebraic K-theory}.

For the proof of assertions
\eqref{the: Induction theorem for algebraic K-theory: calh and Q subset R, R regular} and
\eqref{the: Induction theorem for algebraic K-theory: calh_p and Q subset R, R regular}
use the fact that for regular $R$ with $\bbQ \subseteq R$ the group $G$ satisfies
the (Fibered) Isomorphism Conjecture (in the range $\le N$) for algebraic $K$-theory
with coefficients in $R$ for the family $\calvcyc$ if and only if
$G$ satisfies the Isomorphism Conjecture for algebraic $K$-theory
with coefficients in $R$ for the family $\calfin$. This is proved in \cite[Proposition 2.7]{Lueck-Reich(2004g)}.
Now we can reduce from $\calfin$ to $\calh$ or to $\cale$, provided
that $\bbC\subseteq R$, or after localization at $p$ to $\calh_p$ or to
$\cale_p$, provided that $\bbC \subseteq R$,
by the Transitivity Principle~\ref{the: transitivity principle of bcsfinal}.
\end{proof}

\typeout{--------------------------  Section 5  ------------------------}

\section{Outline of the Proof of the Main Result for Algebraic $L$-Theory}
\label{sec: Outline of the Proof of the Main Result for Algebraic $L$-Theory}

We mention that for $L$-theory we always use the decoration $-\infty $, the Isomorphism
Conjecture is not true for the other decorations such as $s$, $h$ and $p$
(see \cite{Farrell-Jones-Lueck(2002)}).

The proof of Theorem~\ref{the: Induction theorem for algebraic L-theory},
which is the $L$-theory version of Theorem~\ref{the: Induction theorem for algebraic K-theory},
is analogous. The only difference is that one
has to take the involutions into account and replace the functor $\Sw^f$
and the results about it due to Swan  \cite[Corollary 4.2]{Swan(1960a)}  by its $L$-theoretic
version denoted by $GW(H;R)$ and studied by Dress in
\cite[Theorem 2]{Dress(1975)}.

The reduction from $\calfin$ to $\calvcyc$ explained in
\cite[Proposition~2.18]{Lueck-Reich(2004g)} works also in the Fibered
case since the relevant UNil-terms for amalgamated products and
HNN-extensions vanish for all groups (see \cite{Cappell(1974c)}).

\typeout{--------------------------  Section 6  ------------------------}

\section{Proof of the Main Result for Topological  $K$-Theory}
\label{sec: Proof of the Main Result for Topological  $K$-Theory}

The proof of Theorem~\ref{the: Induction theorem for topological K-theory} will
require two ingredients: a Completion Theorem and
a Universal Coefficient Theorem for topological $K$-theory.

In \cite[Theorem 6.5]{Lueck-Oliver(2001b)} a
family version of the Atiyah-Segal Completion Theorem is proved.
In the special case, where $G$ is finite, $X = \pt$ and the family
is the family $\calfcyc$ of finite cyclic subgroups it yields an isomorphism
of pro-groups
\begin{multline}
\{KF_G^*(\pr^{(n)})\}_{n \ge 0} \colon \{KF_G^*(\pt)\}_{n \ge 0}
      \xrightarrow{\cong} \{KF_G^*(\EFG{\calfcyc}{G}^{(n)})\}_{n \ge 0},
\label{Completion result}
\end{multline}
where the first pro-group is given by the constant system with
$ KF_G^*(\pt)$ as value, $\EFG{\calfcyc}{G}^{(n)}$ is the $n$-skeleton of the $G$-$CW$-complex
$\EFG{\calfcyc}{G}$ and the $G$-map $\pr^{(n)} \colon \EFG{\calfcyc}{G}^{(n)} \to G/G$ is the projection.
This means that for every $n$ there is $m > n$ and $f_{m,n} : KF_G^*(\pt) \to KF_G^*(\EFG{\calfcyc}{G}^{(m)})$
such that
\[
\xymatrix
{
KF_G^*(\EFG{\calfcyc}{G}^{(n)}) \ar[d]_{\pr^{(n)}} &
KF_G^*(\EFG{\calfcyc}{G}^{(m)}) \ar[l] \ar[d]^{\pr^{(n)}}
\\
KF_G^*(\pt) &
KF_G^*(\pt) \ar[l]^= \ar[lu]^{f_{m,n}}
}
\]
commutes.
The point is that no $I$-adic completion occurs since the map
$$R_{\bbC}(G) = KF^0_G(\pt) \to \prod_{(C), C \text{ cyclic}} KF_H^0(\pt) = R_{\bbC}(C)$$
is injective and hence $I$ is the zero ideal.
Here $KF$ is either complex $K$-theory ($KF = K$) or real $K$-theory ($KF = KO$).
The complex case occurs already in \cite[Theorem 5.1]{Jackowski(1985)}.

Let $G$ be a finite group.
In \cite{Boekstedt(1981)} B\"okstedt proves Universal Coefficient Theorems
that express equivariant $K$-cohomology in terms of equivariant $K$-homology.
Using S-duality \cite[Chapter~XVI.7]{May(1996)} his results provide also Universal Coefficient Theorems that express
equivariant $K$-homology in terms of equivariant $K$-cohomology.
Let $X$ be a finite $G$-$CW$-complex.
For complex $K$-theory B\"okstedt's result asserts that there are natural short exact sequences
\begin{multline}
0 \to \Ext^1_{K^*_G(\pt)}(K_{*}^G(X), K_*^G(\pt)) \to K^*_G(X) \\
      \to \Hom_{K^*_G(\pt)}(K_*^G(X),K_*^G(\pt)\to 0,
\end{multline}
\vspace{-6ex}
\begin{multline}
\label{eq:UCT-K-Homology}
0 \to \Ext^1_{K^*_G(\pt)}(K^{*}_G(X), K^*_G(\pt)) \to K^G_*(X) \\
      \to \Hom_{K^*_G(\pt)}(K^*_G(X),K^*_G(\pt)\to 0.
\end{multline}
For $KO$-theory his results provide spectral sequences
\begin{eqnarray}
E_2^{p,*} = \Ext^p_{KO^*_G(\pt)}(KO_*^G(X), KO_*^G(\pt)) & \Longrightarrow & KO^*_G(X), \\
\label{eq:UCT-KO-Homology}
E^2_{p,*} = \Ext^p_{KO^*_G(\pt)}(KO^*_G(X), KO_G^*(\pt)) & \Longrightarrow & KO_*^G(X)
\end{eqnarray}
such that $E_\infty^{p,*} = 0$ and $E^\infty_{p,*} = 0$ for $p \geq 2$.
In \cite{Joachim-Lueck(2004)}  the Universal Coefficient Theorem for complex $K$-theory is generalized to
infinite groups and proper $G$-$CW$-complexes.

We can now give the proof of Theorem~\ref{the: Induction theorem for topological K-theory}.

\begin{proof}
Because of the Transitivity Principle~\ref{the: transitivity} it suffices
to show that
\begin{eqnarray}
KF^G_n(\pr) \colon KF^G_n(\EFG{\calfcyc}{G}) &\xrightarrow{\cong}& KF^G_n(\pt)
\label{topological K-theory and calfcyc}
\end{eqnarray}
is bijective for all finite groups $G$ and $n \in \bbZ$.
There exists a $G$-$CW$-model for $\EFG{\calfcyc}{G}$ whose skeleta are all finite $G$-$CW$-complexes.
This follows for instance from the functorial construction in
\cite[Section 3 and Lemma 7.6]{Davis-Lueck(1998)} using the fact that $\OrF{\calf}{G}$ is a category
with finitely many morphisms.

If we apply \eqref{eq:UCT-K-Homology} to $\pr^{(n)} \colon \EFG{\calfcyc}{G}^{(n)} \to G/G$
we obtain a map between two short exact sequences.
Since $\colim_{n \to \infty}$ is an exact functor, these sequences stay exact if we apply $\colim_{n \to \infty}$.
Because of the isomorphism of pro-groups \eqref{Completion result},
we get isomorphism in the first and third term of \eqref{eq:UCT-K-Homology}
(see for instance \cite{Joachim-Lueck(2004)}).
By the Five-Lemma this implies that
\[
\colim_{n \to \infty} K_*^G(\EFG{\calfcyc}{G}^{(n-1)}) \to \colim_{n \to \infty}   K^G_*(\pt)
\]
is bijective. Since $K$-homology is compatible
with colimits, this map can be identified with the map
\eqref{topological K-theory and calfcyc} for complex $K$-theory.

For $KO$-theory we can proceed similarly.
The isomorphism of pro-groups \eqref{Completion result} yields
isomorphisms on the colimit of the $E^2$-Term of the spectral sequence \eqref{eq:UCT-KO-Homology}.
Since $\colim_{n \to \infty}$ is an exact functor this yields also isomorphisms on the colimit of the $E^\infty$ term.
Because this $E^\infty$ term has only a finite number of lines a simple diagram chase shows that the
induced map
\[
\colim_{n \to \infty} KO_*^G(\EFG{\calfcyc}{G}^{(n)})  \to \colim_{n \to \infty}   KO^G_*(\pt)
\]
on the right hand side of \eqref{eq:UCT-KO-Homology} is an isomorphisms.
Since $KO$-homology is compatible
with colimits, this map can be identified with the map
\eqref{topological K-theory and calfcyc} for $KO$-theory.
\end{proof}

\typeout{--------------------------  Section 7  ------------------------}

\section{Versions in Terms of Colimits}
\label{sec: Versions in terms of colimits}

The classical induction theorems can be stated in terms of colimits.
If we combine Theorem~\ref{the: Criterion for induction for equivariant homology theories}
and Lemma~\ref{lem: calh^?_*(-,bfK) and assumptions for induction} and (its $L$-theory version)
we get for an extension $1 \to K \to G \xrightarrow{p} F \to 1$ for a finite group $F$ isomorphisms
\begin{eqnarray*}
\colim_{F/L \in \OrF{\calh}{F}} K_n(R[p^{-1}(L)]) & \xrightarrow{\cong} & K_n(RG);
\\
\colim_{F/L \in \OrF{\calh_p}{F}} K_n(R[p^{-1}(L)])_{(p)} & \xrightarrow{\cong} & K_n(RG)_{(p)};
\\
\colim_{F/L \in \OrF{\calh_2 \cup \bigcup_{p\text{ prime}, p \not = 2} \cale_p}{F}}
L^{\langle - \infty \rangle}_n(R[p^{-1}(L)]) & \xrightarrow{\cong} & L_n^{\langle - \infty \rangle}(RG);
\\
\colim_{F/L \in \OrF{\bigcup_{p\text{ prime}, p \not = 2} \cale_p}{F}} L_n(R[p^{-1}(L)])[1/2] &
         \xrightarrow{\cong} & L_n(RG)[1/2].
\end{eqnarray*}
But in general our results cannot be stated in such an elementary way. For instance,
Theorem~\ref{the: Induction theorem for topological K-theory} says for a finite group $G$
that there is an isomorphism
$$K_0^G(\EFG{\calfcyc}{G}) \xrightarrow{\cong} K_0^G(G/G)$$
In general $K_0^G(G/H)$ is the complex representation ring $R_{\bbC}(H)$ of the finite group $H$.
It is not true in general that the canonical map
$$\colim_{G/H \in \OrF{\calfcyc}{G}} R_{\bbC}(H) ~ \to ~ R_{\bbC}(G)$$
is bijective or, equivalently, that the canonical map
$$\colim_{G/H \in \OrF{\calfcyc}{G}} K_0^G(G/H) ~ \to ~ K_0^G(G/G)$$
is bijective.

In some special cases one can get formulations of our results in terms of colimits.

\begin{theorem}
\label{the: special results for K_1 and K_0}
\begin{enumerate}
\item \label{the: special results for K_1 and K_0: K_{-1}}
The group $G$ satisfies the Isomorphism Conjecture in the range $n \le -1$ for algebraic $K$-theory
with coefficients in $R = \bbZ$ for the family $\calvcyc$ if and only if
$$K_n(\bbZ G ) = 0 \text{ for  } n \le -2$$
and the canonical map
$$\colim_{G/H \in \OrF{\calh}{G}} K_{-1}(\bbZ H) ~ \xrightarrow{\cong} ~ K_{-1}(\bbZ G)$$
is bijective;

\item \label{the: special results for K_1 and K_0: K_0}
Suppose that $R$ is regular and $\bbQ \subseteq R$. Then the  group $G$ satisfies the Isomorphism
Conjecture in the range $n \le 0$ for algebraic $K$-theory
with coefficients in $R$ for the family $\calvcyc$ if and only if
$$K_n(RG ) = 0 \text{ for  } n \le -1$$
and the canonical map
$$\colim_{G/H \in \OrF{\calh}{G}} K_{0}(RH) ~ \xrightarrow{\cong} ~ K_0(RG)$$
is bijective.

If we assume that $R$ is regular and $\bbC \subseteq R$, then we can
replace $\calh$ by $\cale$.
\end{enumerate}
\end{theorem}

\begin{proof}
\eqref{the: special results for K_1 and K_0: K_{-1}}
It follows from \cite{Farrell-Jones(1995)} that for a virtually cyclic group $V$ the assembly map
$$ H_n^G(\EFG{\calfin}{V};\bfK_{\bbZ}) ~ \xrightarrow ~ K_n(\bbZ V)$$
is surjective in the range $n \leq -1$.
Moreover, this map is  known to be injective in all degrees \cite{Bartels(2003h)}, \cite{Rosenthal(2003)}.
Theorem~\ref{the: Criterion for induction for equivariant homology theories} and
Lemma~\ref{lem: calh^?_*(-,bfK) and assumptions for induction} imply that
for every finite group $F$ the assembly map
$$H_n^F(\EFG{\calh}{F};\bfK_{\bbZ}) ~ \xrightarrow{\cong} K_n(\bbZ F)$$
is bijective.
Applying the Transitivity Principle for equivariant homology~\ref{the: transitivity principle of bcsfinal} twice
we conclude that a group $G$ satisfies the Isomorphism
Conjecture for algebraic $K$-theory
with coefficients in $R = \bbZ $ for the family $\calvcyc$ in the range $n \le -1$
if and only if  the assembly map
$$H_n^G(\EFG{\calh}{G};\bfK_{\bbZ}) ~ \xrightarrow{\cong} ~ K_n(\bbZ G)$$
is bijective for $n \le -1$.
It is proven in \cite{Carter(1980)} that for finite groups $F$ and $n \leq -2$ $K_n(\bbZ F) = 0$.
Now  a spectral sequence argument shows
that $H_n^G(\EFG{\calh}{G};\bfK_{\bbZ})$ vanishes for $n \le -2$ and can be identified
with $\colim_{G/H \in \OrF{\calh}{F}} K_{-1}(\bbZ H)$ for $n = -1$.
\\[1mm]
\eqref{the: special results for K_1 and K_0: K_0}
By Theorem~\ref{the: Induction theorem for algebraic K-theory}
\eqref{the: Induction theorem for algebraic K-theory: calh and Q subset R, R regular}
the  group $G$ satisfies the Isomorphism
Conjecture in the range $n \le 0$ for algebraic $K$-theory
with coefficients in $R$ for the family $\calvcyc$ in the range $n \le 0$ if and only if
the canonical map
$$H_n^G(\EFG{\calh}{G};\bfK_R) ~ \xrightarrow{\cong} ~ K_n(RG)$$
is bijective for $n \le 0$. Since $RH$ is regular for a finite group $H$,
$K_q(RH) = 0$ for $q \le -1$. Now  a spectral sequence argument shows
that $H_n^G(\EFG{\calh}{G};\bfK_R)$ vanishes for $n \le -1$ and can be identified
with $\colim_{G/H \in \OrF{\calh}{G}} K_{0}(RH)$ for $n = 0$.
\end{proof}

\typeout{--------------------------  Section 8 ------------------------}

\section{On Quinn's Hyperelementary Induction Conjecture}
\label{sec: On Quinn's Hyperelementary Induction Conjecture}

Quinn \cite{Quinn(2004)} states the following conjecture

\begin{conjecture}[Hyperelementary Induction Conjecture]
\label{con: Hyperelementary induction conjecture}
All groups satisfy induction for the class $\overline{\calh}$ of (possibly infinite) hyperelementary groups.
\end{conjecture}

The phrase that a group $G$ satisfies induction for a family $\calf$ of subgroups of $G$ means
that the Fibered Farrell-Jones Conjecture is true for the family $\calf$
and algebraic $K$-theory with coefficients in
a given associative ring with unit $R$. The Fibered Isomorphisms Conjecture is due to
Farrell-Jones~\cite{Farrell-Jones(1993a)} (see also \cite[Section 4.2.2]{Lueck-Reich(2004g)}).
Quinn calls a group $H$ \emph{$p$-hyperelementary}
if there exists a prime $p$, a finite $p$-group $P$ and a (possibly infinite) cyclic group $C$
such that $H$ can be written as an extension $1 \to C \to H \to P \to 1$, and hyperelementary
if it is $p$-hyperelementary for some prime $p$. The class $\overline{\calh}$ appearing
in the Hyperelementary Induction Conjecture~\ref{con: Hyperelementary induction conjecture}
is to be understood with respect to Quinn's definition of hyperelementary group as above.
Of course $\overline{\calh} \cap \calfin$ is the class $\calh$ of finite hyperelementary groups.
There is an inclusion $\overline{\calh} \subseteq \calh'$ but the classes
$\overline{\calh}$ and $\calh'$ are different. One would like to prove
that Hyperelementary Induction Conjecture~\ref{con: Hyperelementary induction conjecture}
follows from the Fibered Isomorphism Conjecture for algebraic $K$-theory. This is in view of
Theorem~\ref{the: Induction theorem for algebraic K-theory}~\eqref{the: Induction theorem for algebraic K-theory: calh}
and Theorem~\ref{the: transitivity} equivalent to a positive answer to the following question:

\begin{question}
\label{que: relating Quinn's Conjecture to our results}
Is for every group $H \in \calh'$ the Hyperelementary
Induction Conjecture~\ref{con: Hyperelementary induction conjecture} true?
\end{question}

We have no evidence for a positive answer from general machinery, a proof of a positive answer will
need some special information about  the Nil-groups $NK_n(RH)$.

\def\cprime{$'$} \def\polhk#1{\setbox0=\hbox{#1}{\ooalign{\hidewidth
  \lower1.5ex\hbox{`}\hidewidth\crcr\unhbox0}}}

\end{document}